\documentclass[10pt]{amsart}
\usepackage{graphicx,amssymb}
\usepackage{a4wide}

\def\baselinestretch{1.2}

\theoremstyle{plain}
\newtheorem{theorem}{Theorem}[section]
\newtheorem{lemma}[theorem]{Lemma}
\newtheorem{corollary}[theorem]{Corollary}
\newtheorem{proposition}[theorem]{Proposition}

\theoremstyle{definition}
\newtheorem{definition}[theorem]{Definition}
\newtheorem{example}[theorem]{Example}

\newtheorem{remark}[theorem]{Remark}

\newenvironment{algorithm}[1]
{\par\medskip\noindent\textbf{#1.}\ \ \ignorespaces}
{\par\medskip}

\def\B{\mathbf B}
\def\D{\mathcal D}
\def\Z{\mathbf Z}

\def\Sym{\operatorname{\mathit{Sym}}}

\def\len{\operatorname{len}}
\def\c{\mathbf c}
\def\d{\mathbf d}

\def\C{\mathcal C}

\def\myangle#1{\langle #1\rangle}
\def\myceil#1{\lceil #1\rceil}
\def\myfloor#1{\lfloor #1\rfloor}

\def\t{\operatorname{\it t}}
\def\INF{\t_{\inf}}
\def\SUP{\t_{\sup}}
\def\LEN{\t_{\len}}
\def\infs{\inf{\!}_ s}
\def\sups{\sup{\!}_s}
\def\lens{\len{\!}_s}

\def\wedgeL{\wedge_L}
\def\wedgeR{\wedge_R}
\def\veeL{\vee_{\!L}}
\def\veeR{\vee_{\!R}}

\def\Artinbr{{\operatorname{[Artin]}}}
\def\BKLbr{{\operatorname{[BKL]}}}

\def\BKL{{\operatorname{BKL}}}
\def\atom{{\operatorname{atom}}}
\def\simple{{\operatorname{simple}}}
\def\Tlat{T_{\operatorname{lattice}}}
\let\epsilon\varepsilon

\begin{document}

\title{Conjugacy classes of periodic braids}

\author{Eon-Kyung Lee and Sang Jin Lee}
\address{Department of Applied Mathematics, Sejong University,
    Seoul, 143-747, Korea}
\email{eonkyung@sejong.ac.kr}
\address{Department of Mathematics, Konkuk University,
    Seoul, 143-701, Korea}
\email{sangjin@konkuk.ac.kr}
\date{\today}

\begin{abstract}
Recently, there have been several progresses
for the conjugacy search problem (CSP) in Garside groups,
especially in braid groups.
All known algorithms for solving this problem use
a sort of exhaustive search in a particular finite set
such as the super summit set and the ultra summit set.
Their complexities are proportional to the size of the finite set,
even when there exist very short conjugating elements.
However, ultra summit sets are very large in some cases
especially for reducible braids and periodic braids.
Some possible approaches to resolve this difficulty would be
either to use different Garside structures and Garside groups
in order to get a sufficiently small ultra summit set,
or to develop an algorithm for finding a conjugating element
faster than exhaustive search.
Using the former method,
Birman, Gonz\'alez-Meneses and Gebhardt have proposed
a polynomial-time algorithm for the CSP for periodic braids.

In this paper we study the conjugacy classes of periodic braids
under the BKL Garside structure,
and show that we can solve the CSP for periodic braids in polynomial time
although their ultra summit sets are exponentially large.
Our algorithm describes how to connect two periodic braids
in the (possibly exponentially large) ultra summit set
by applying partial cycling polynomially many times.

\end{abstract}

\maketitle

\section{Introduction}

The conjugacy problem in groups has two versions:
the conjugacy decision problem (CDP) is to decide
whether given two elements are conjugate or not;
the conjugacy search problem (CSP) is to find a conjugating element
for a given pair of conjugate elements.
The conjugacy problem is of great interest for the Artin braid group $B_n$,
which has the well-known Artin presentation~\cite{Art47}:
$$
B_n  =  \left\langle \sigma_1,\ldots,\sigma_{n-1} \left|
\begin{array}{ll}
\sigma_i \sigma_j = \sigma_j \sigma_i & \mbox{if } |i-j| > 1, \\
\sigma_i \sigma_j \sigma_i = \sigma_j \sigma_i \sigma_j & \mbox{if } |i-j| = 1.
\end{array}
\right.\right\rangle.
$$

In the late sixties, Garside~\cite{Gar69} first solved the conjugacy problem
in the braid groups, and then his theory has been generalized
and enriched by many mathematicians.
The algorithms for solving the conjugacy problem, provided by Garside's theory,
involve computation of finite nonempty subsets of the conjugacy class
such as the summit set~\cite{Gar69},
the super summit set~\cite{EM94, BKL98, FG03},
the ultra summit set~\cite{Geb05},
the stable super summit set~\cite{LL06c},
the stable ultra summit set~\cite{BGG06b}
and the reduced super summit set~\cite{KL06}.
Let us call these sets \emph{conjugacy representative sets}.
The conjugacy representative sets depend not only on the group itself
but also on the Garside structure on it,
which is a pair of a positive monoid
(equivalently, a lattice order invariant under left multiplication)
and a Garside element.
The braid group $B_n$ admits two well-known
Garside structures, the Artin Garside structure
and the BKL Garside structure.
Let $B_n^\BKLbr$ and $B_n^\Artinbr$ denote the braid group $B_n$
endowed with the Artin and the BKL Garside structure, respectively.

The best known upper bound of the complexity for computing a conjugacy
representative set $S$ is of the form $\mathcal O(|S|\cdot p(n))$,
where $|S|$ denotes the size of the set $S$ and $p(n)$ is a polynomial in $n$.
Franco and Gonz\'alez-Meneses~\cite{FG03} and Gebhardt~\cite{Geb05}
showed this when $S$ is the super summit set
and the ultra summit set, respectively,
and it follows easily from their results that
the complexities for computing the other conjugacy representative
sets have upper bound of the same form.

All of the known algorithms for solving the CSP in braid groups,
and more generally in Garside groups, use a sort of \emph{exhaustive search}
in conjugacy representative sets.
The most popular method is as follows:
given two elements $\alpha$ and $\beta$, one computes
the conjugacy representative set $S(\alpha)$ of $\alpha$
and an element $\beta'$ in the conjugacy representative set of $\beta$,
and then check whether $\beta'$ belongs to $S(\alpha)$.
Therefore, the complexity of this kind of algorithms for the CSP
is at least the complexity for computing conjugacy representative sets,
even when there exists a very short conjugating element.
Recently, Birman, Gebhardt and Gonz\'alez-Meneses~\cite{BGG06b}
introduced two classes of subsets, called black and grey components,
of an ultra summit set, and showed that the conjugacy problem can be solved
by computing one black component and one grey component.
This new algorithm is faster than the previous ones,
however it does not improve the theoretical complexity.

The size of a conjugacy representative set is exponential
in the braid index $n$ in some cases,
especially for reducible and periodic braids.
See~\cite[\S2]{BGG06c} or Remark~\ref{rmk:size}.
(An $n$-braid is said to be {\em periodic}
if some power of it belongs to $\langle\Delta^2\rangle$, the center of $B_n$,
and \emph{reducible} if there is an essential curve system in the
punctured disk which is invariant under the action of the braid.)
Therefore, a possible way to solve the CSP more efficiently
would be either to use different Garside structures and Garside groups
such that the conjugacy representative set in question is small enough,
or to develop an algorithm for computing a conjugating element
not in an exhaustive way but in a deterministic way.

For periodic braids, Birman, Gonz\'alez-Meneses and Gebhardt showed
in~\cite{BGG06b} that
their new algorithm in~\cite{BGG06b} is not better than the previous
one in~\cite{Geb05}.
Hence they used in~\cite{BGG06c} different Garside structures and Garside groups
in order to get a small super summit set.
Using this method, they have proposed
a polynomial-time algorithm for the CSP for periodic braids.

It is well-known that the CDP for periodic braids is very easy:
the results of Eilenberg~\cite{Eil34} and Ker\'ekj\'art\'o~\cite{Ker19}
imply that an $n$-braid is periodic if and only if
it is conjugate to a power of either $\delta$ or $\epsilon$,
where $\delta =\delta_{(n)} = \sigma_{n-1}\sigma_{n-2}\cdots\sigma_1$ and
$\epsilon = \epsilon_{(n)}=\delta\sigma_1$, hence
\begin{itemize}
\item an $n$-braid $\alpha$ is periodic if and only if
either $\alpha^n$ or $\alpha^{n-1}$ belongs to $\langle\Delta^2\rangle$;
\item two periodic braids are conjugate if and only if
they have the same exponent sum.
\end{itemize}
In~\cite{BGG06c}, Birman, Gonz\'alez-Meneses and Gebhardt showed the following.

\begin{itemize}
\item
The sizes of the ultra summit sets of $\delta$ and $\epsilon$
in $B_n^\Artinbr$  are exponential in the braid index $n$,
hence the complexity of usual algorithms for the CSP is not
polynomial for periodic braids.

\item
The super summit set of $\delta^k$ in $B_n^\BKLbr$ is $\{\delta^k\}$,
hence the CSP for $\delta$-type periodic $n$-braids
is solvable in polynomial time (in the braid index $n$ and the input word length)
by using the BKL Garside structure on $B_n$.

\item
The CSP for $\epsilon$-type periodic $n$-braids
is solvable in polynomial time (in the braid index $n$ and the input word length)
by using an algorithm for the CSP
for $\delta$-type periodic braids in $B_{2n-2}^\BKLbr$,
together with algorithms for computing an isomorphism from
$P_{n,2}$ to $\Sym_{2n-2}$ and its inverse,
where $P_{n,2}$ is the subgroup of $B_n$ consisting of $n$-braids
that fix the second puncture and $\Sym_{2n-2}$ is the centralizer
of $\delta_{(2n-2)}^{n-1}$ in $B_{2n-2}$.
\end{itemize}

In this paper we develop a new polynomial-time
(in the braid index and the input word length) algorithm for solving
the CSP for periodic braids only by exploiting the BKL Garside structure
on $B_n$, and study how to improve efficiency of the algorithms.
Compared to the algorithm of Birman et.~al.~in~\cite{BGG06c},
our algorithm has lower complexity and, moreover,
it uses a single Garside group $B_n^\BKLbr$, hence the implementation
is simpler.

First, we study periodic elements in Garside groups.
An element of a Garside group is said to be \emph{periodic} if
some power of it belongs to the cyclic group generated by the Garside
element.
For periodic elements in Garside groups, the super summit set
is the same as the ultra summit set.
Observing the previous results of Bestvina~\cite{Bes99} and
Charney, Meier and Whittlesey~\cite{CMW04},
we can see that every periodic element in Garside groups has
a special type of power, which we will call a \emph{BCMW-power}.
We show the characteristics of BCMW-powers and its interesting property:
the CSP for two periodic elements $g$ and $h$ in a Garside group
is equivalent to the CSP for $g^r$ and $h^r$, where $g^r$ is a
BCMW-power of $g$.
Then we study super summit sets of periodic elements in Garside groups.
Especially, we show that super summit sets of a certain type of
periodic elements are closed under any partial cycling.

We show that, for some integers $k$,
the super summit set of $\epsilon^k$ in $B_n^\BKLbr$
is exponentially large in the braid index $n$.
Using the results on periodic elements in Garside groups and
using the characteristics of the BKL Garside structure on braid groups,
we present an explicit method to transform an arbitrary braid in
the super summit set of $\epsilon^k$ to $\epsilon^k$
by applying partial cycling polynomially many times,
even when the super summit set is exponentially large.

On the other hand, we discuss concrete methods for improving efficiency
of the algorithms.
Especially, using a known algorithm for powering integers~\cite{Coh93,Sho05}
and our recent result on
abelian subgroups of Garside groups~\cite{LL06c}, we propose an algorithm
for computing a super summit element of a power of a periodic element
in a Garside group more efficiently.

We hope that the results of this paper on periodic elements in Garside groups
are also useful in studying periodic elements in other Garside groups.
This paper is organized as follows.
\S2 gives a brief review on Garside groups and
braids groups with the BKL Garside structure.
\S3 studies periodic elements in Garside groups and their super
summit sets.
\S4 studies super summit sets of $\epsilon$-type periodic braids in $B_n^\BKLbr$,
and shows a method to find a conjugating element for any two braids
in the super summit set of $\epsilon^k$.
\S5 constructs algorithms for the conjugacy problem for periodic braids
in $B_n^\BKLbr$.

\section{Preliminaries}

\subsection{Garside monoids and groups}

The class of Garside groups, first introduced by Dehornoy and Paris~\cite{DP99},
provides a lattice-theoretic generalization of the braid groups
and the Artin groups of finite type.

For a monoid $M$, let $e$ denote the identity element.
An element $a\in M\setminus \{e\}$ is called an \emph{atom} if
$a=bc$ for $b,c\in M$ implies either $b=e$ or $c=e$.
For $a\in M$, let $\Vert a\Vert$ be the supremum
of the lengths of all expressions of
$a$ in terms of atoms. The monoid $M$ is said to be \emph{atomic}
if it is generated by its atoms and $\Vert a\Vert<\infty$
for any element $a$ of $M$.
In an atomic monoid $M$, there are partial orders $\le_L$ and $\le_R$:
$a\le_L b$ if $ac=b$ for some $c\in M$;
$a\le_R b$ if $ca=b$ for some $c\in M$.

\begin{definition}
An atomic monoid $M$ is called a \emph{Garside monoid} if
\begin{enumerate}
\item[(i)] $M$ is finitely generated;
\item[(ii)] $M$ is left and right cancellative;
\item[(iii)] $(M,\le_L)$ and $(M,\le_R)$ are lattices;
\item[(iv)] $M$ contains an element $\Delta$, called a
\emph{Garside element}, satisfying the following:\\
(a) for each $a\in M$, $a\le_L\Delta$ if and only if $a\le_R\Delta$;\\
(b) the set $\{a\in M: a \le_L\Delta\}$ generates $M$.
\end{enumerate}
\end{definition}

Recall that a partially ordered set $(P,\le)$ is called a
\emph{lattice} if there exist the gcd $a\wedge b$ and the lcm $a\vee b$
for any $a,b\in P$.
The gcd $a\wedge b$ is the unique element such that
(i) $a\wedge b\le a$ and $a\wedge b\le b$;
(ii) if $c$ is an element satisfying $c\le a$ and $c\le b$,
then $c\le a\wedge b$.
Similarly, the lcm $a\vee b$ is the unique element such that
(i) $a\le a\vee b$ and $b\le a\vee b$;
(ii) if $c$ is an element satisfying $a\le c$ and $b\le c$,
then $a\vee b\le c$.
Let $\wedgeL$ and $\veeL$ (resp. $\wedgeR$ and $\veeR$) denote
the gcd and lcm with respect to $\le_L$ (resp. $\le_R$).

An element $a$ of $M$ is called a \emph{simple element} if $a\le_L\Delta$.
Let $\D$ denote the set of all simple elements.
A \emph{Garside group} is defined as the group of fractions
of a Garside monoid.
When $M$ is a Garside monoid and $G$ is the group of fractions of $M$,
we identify the elements of $M$ and their images in $G$
and call them \emph{positive elements} of $G$.
$M$ is called the \emph{positive monoid} of $G$,
often denoted $G^+$.
The triple $(G, G^+, \Delta)$ is called a
\emph{Garside structure} on $G$.
We remark that a given group $G$ may
admit more than one Garside structures.

The partial orders $\le_L$ and $\le_R$, and thus the lattice structures
in the positive monoid $G^+$ can be extended to the Garside group $G$.
For $g, h\in G$, $g\le_L h$ (resp. $g\le_R h$) means $g^{-1}h\in G^+$
(resp. $hg^{-1}\in G^+$), in which
$g$ is called a {\em prefix} (resp. {\em suffix}) of $h$.

For an element $g\in G$, the \emph{(left) normal form} of $g$ is
the unique expression
$$
g=\Delta^u a_1\cdots a_\ell,
$$
where $u=\max\{r\in\Z:\Delta^r\le_L g\}$,
$a_1,\ldots,a_\ell\in \D\setminus\{e,\Delta\}$ and
$(a_i a_{i+1}\cdots a_\ell)\wedgeL \Delta=a_i$ for $i=1,\ldots,\ell$.
In this case, $\inf(g)=u$, $\sup(g)=u+\ell$ and \/$\len(g)=\ell$ are
called the {\em infimum}, {\em suprimum} and {\em canonical length} of $g$,
respectively.

Let $\tau : G\to G$ be the inner automorphism of $G$
defined by $\tau(g)=\Delta^{-1}g\Delta$ for all $g\in G$.
The {\em cycling} and {\em decycling} of $g$ are defined as
\begin{eqnarray*}
\c(g)&=& \Delta^u a_2\cdots a_\ell\tau^{-u}(a_1),\\
\d(g)&=&\Delta^u \tau^{u}(a_\ell)a_1\cdots a_{\ell-1}.
\end{eqnarray*}

We denote the conjugacy class $\{ h^{-1}gh : h\in G\}$ of $g$ in $G$
by $[g]$.
Define $\infs(g)=\max\{\inf(h):h\in [g]\}$,
$\sups(g)=\min\{\sup(h):h\in [g]\}$ and $\lens(g)=\sups(g)-\infs(g)$,
which are called the {\em summit infimum}, {\em summit suprimum} and
{\em summit length} of $g$, respectively.
The \emph{super summit set} $[g]^S$,
the \emph{ultra summit set} $[g]^U$
and the \emph{stable super summit set} $[g]^{St}$
are finite, nonempty subsets of the conjugacy class of $g$
defined as follows (see~\cite{FG03,Geb05,LL06c} for more detail):
\begin{eqnarray*}
[g]^S&=&\{h\in [g]:\inf(h)=\infs(g) \ \mbox{ and } \sup(h)=\sups(g)\};\\{}
[g]^{U}&=&\{h\in [g]^S: \c^k(h)=h \ \mbox{ for some positive integer $k$}\};\\{}
[g]^{St}&=&\{h\in [g]^S:h^k\in[g^k]^S \ \mbox{ for all positive integers $k$}\}.
\end{eqnarray*}
Elements of super summit sets are called \emph{super summit elements}.

For every element $g\in G$, the following limits are well-defined:
$$
\INF(g)=\lim_{n\to\infty}\frac{\inf(g^n)}n;\quad
\SUP(g)=\lim_{n\to\infty}\frac{\sup(g^n)}n;\quad
\LEN(g)=\lim_{n\to\infty}\frac{\len(g^n)}n.
$$
These limits were introduced in~\cite{LL06a}
to investigate translation numbers in Garside groups.
We will exploit the following properties especially in \S\ref{sec:partial_decycling}.

\begin{proposition}[{\cite{LL06a,LL06b}}]
For $g, h\in G$,
\begin{enumerate}
\item[(i)]
$\INF(h^{-1}gh)=\INF(g)$ and $\SUP(h^{-1}gh)=\SUP(g)$;

\item[(ii)]
$\INF(g^n)= n\cdot\INF(g)$ and $\SUP(g^n)= n\cdot\SUP(g)$
for all integers $n\ge 1$;

\item[(iii)]
$\infs(g)=\lfloor \INF(g)\rfloor$ and $\sups(g)=\lceil \SUP(g)\rceil$;

\item[(iv)]
$\INF(g)$ and $\SUP(g)$ are rational of the form $p/q$,
where $p$ and $q$ are relatively prime integers and
$1\le q\le\Vert\Delta\Vert$.
\end{enumerate}
\end{proposition}

\subsection{The BKL Garisde structure on braid groups}

Birman, Ko and Lee~\cite{BKL98} introduced
a then-new monoid---with the following explicit presentation---whose
group of fractions is the braid group $B_n$:
$$
B_n  =  \left\langle a_{ij}, \ 1\le j < i\le n \left|
\begin{array}{ll}
a_{kl}a_{ij}=a_{ij}a_{kl} & \mbox{if $(k-i)(k-j)(l-i)(l-j)>0$}, \\
a_{ij}a_{jk}=a_{jk}a_{ik}=a_{ik}a_{ij} & \mbox{if $1\le k<j<i \le n$}.
\end{array}
\right.\right\rangle.
$$
The generators $a_{ij}$ are called {\em band generators}.
They are related to the classical generators by
$a_{ij}=\sigma_{i-1}\sigma_{i-2}\cdots\sigma_{j+1}\sigma_j
\sigma_{j+1}^{-1}\cdots\sigma_{i-2}^{-1}\sigma_{i-1}^{-1}$.
The BKL Garside structure of $B_n$ is determined by
the positive monoid which consists of the elements
represented by positive words in the band generators
and the Garside element
$$
\delta = a_{n,n-1}\cdots a_{3,2} a_{2,1}.
$$

\begin{figure}
\includegraphics[scale=1.1]{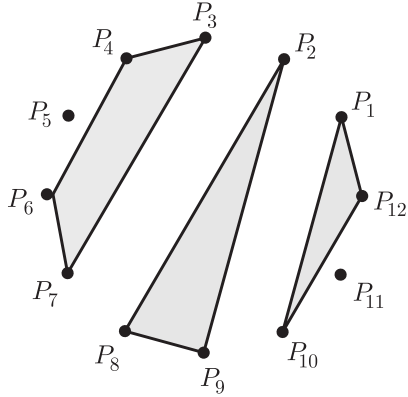}
\caption{The shaded regions show the blocks in the non-crossing
partition corresponding to the simple element
$[12,10,1]\,[9,8,2]\,[7,6,4,3]$ in $B_{12}$.}\label{fig:ncp}
\end{figure}

The simple elements in the BKL Garside structure are in one-to-one
correspondence with non-crossing partitions.
(Note that the simple elements in $B_n^\Artinbr$ are
in one-to-one correspondence with $n$-permutations.)
Let $P_1,\ldots,P_n$ be the points in
the complex plain given by $P_k=\exp(\frac{2k\pi}ni)$.
See Figure~\ref{fig:ncp}.
Recall that a partition of a set is a collection of pairwise
disjoint subsets whose union is the entire set.
Those subsets (in the collection) are called {\em blocks}.
A partition of $\{P_1,\ldots,P_n\}$
is called a \emph{non-crossing partition}
if the convex hulls of the respective blocks are pairwise disjoint.

A positive word of the form
$a_{i_k i_{k-1}} \cdots a_{i_3 i_2} a_{i_{2}i_1}$,
$1\le i_1 < i_2 < \cdots < i_k\le n$, is called a \emph{descending cycle}
and denoted $[i_k,\ldots,i_2,i_1]$.
Two descending cycles, $[i_k,\ldots,i_1]$ and $[j_l,\ldots,j_1]$,
are said to be \emph{parallel}
if the convex hulls of $\{P_{i_1},\ldots,P_{i_k}\}$
and of $\{P_{j_1},\ldots,P_{j_l}\}$ are disjoint.
A simple element is a product of parallel descending cycles.

\section{Super summit sets of periodic elements in Garside groups}
\label{sec:partial_decycling}

In this section $G$ denotes a Garside group with Garside element $\Delta$.
Let $m$ denote the smallest positive integer such that
$\Delta^m$ is central in $G$,
and let $G_\Delta$ denote the central quotient $G/\myangle{\Delta^m}$
where $\myangle{\Delta^m}$ is the cyclic group generated by $\Delta^m$.
For an element $g\in G$, let $\bar g$ denote the image of $g$
under the natural projection from $G$ to $G_\Delta$.
An element $g$ of $G$ is said to be \emph{periodic} if
$\bar g$ has a finite order in $G_\Delta$ or, equivalently,
if some power of $g$ belongs to $\langle\Delta^m\rangle$.
In this section we study periodic elements in Garside groups.

\subsection{BCMW powers of periodic elements}
Due to the work of Dehornoy~\cite{Deh98}, we know that
Garside groups are torsion free,
hence there is no non-trivial finite subgroup.
We start this section with the results of Bestvina~\cite{Bes99} and
Charney, Meier and Whittlesey~\cite{CMW04}
for finite subgroups of $G_\Delta$.
The following theorem was first proved
by Bestvina~\cite[Theorem 4.5]{Bes99} for
Artin groups of finite type and then generalized to Garside groups by
Charney, Meier and Whittlesey~\cite[Corollary 6.8]{CMW04}.

\begin{theorem}
[Bestvina~\cite{Bes99}, Charney, Meier and Whittlesey~\cite{CMW04}]
\label{thm:CMW}
The finite subgroups of\/ $G_\Delta$ are, up to conjugacy,
one of the following two types:
\begin{itemize}
\item[(i)]
the cyclic group generated by $\Delta^u a$ for an integer $u$ and
a simple element $a\in\D\backslash\{\Delta\}$
such that if $a\neq e$, then for some integer $q$
with $2\le q\le\Vert\Delta\Vert$,
\begin{equation}\label{eq:CMW}
\tau^{(q-1)u}(a)\,\tau^{(q-2)u}(a)\cdots \tau^{u}(a)\,a=\Delta ;
\end{equation}

\item[(ii)]
the direct product of a cyclic group of type (i) and
$\langle\Delta^k\rangle$ where $\Delta^k$ commutes with $a$.
\end{itemize}
\end{theorem}

Notice that the element $\Delta^u a$ in Theorem~\ref{thm:CMW}~(i)
is clearly a periodic element because
$(\Delta^u a)^q
=\Delta^{qu}\,\tau^{(q-1)u}(a)\cdots \tau^{u}(a)\,a
=\Delta^{qu+1}$.
However, not every periodic element satisfies the conditions
in Theorem~\ref{thm:CMW}~(i).
Namely, if $\Delta^ua $ is periodic, then
$\tau^{(q-1)u}(a)\,\tau^{(q-2)u}(a)\cdots \tau^{u}(a)\,a=\Delta^t$
for some positive integers $q$ and $t$, but $t$ is not necessarily 1.
Motivated by this observation,
we define the following notions for periodic elements.

\begin{definition}\label{def:BCMW}
Let $g$ be a periodic element of a Garside group $G$ such that $\INF(g)=p/q$
for relatively prime integers $p$ and $q$ with $q\ge 1$.
\begin{itemize}
\item[(i)]
The periodic element $g$ is said to be \emph{P-minimal} if $p\equiv 1\bmod q$.
\item[(ii)]
A power $h=g^r$ is called a \emph{BCMW-power} of $g$
if $h$ is P-minimal and
$\bar h$ generates the same cyclic subgroup of $G_\Delta$ as $\bar g$ does.
\item[(iii)]
The periodic element $g$ is said to be \emph{C-tight}
if $p\equiv 0\bmod m$.
\end{itemize}
\end{definition}

It will be shown in Lemma~\ref{lem:Primi_equiv}
that a periodic element $g$ is P-minimal
if and only if $g$ is conjugate to an element of the form $\Delta^u a$
satisfying the conditions in Theorem~\ref{thm:CMW}~(i).
(It will be also shown that if $\Delta^u a$ is P-minimal and belongs to
its stable super summit set, then the simple element $a$ is minimal
among the positive parts of the powers of $\Delta^u a$ other than
the identity.)
Therefore, if restricted to finite cyclic subgroups of
$G_\Delta$, Theorem~\ref{thm:CMW} is equivalent to the existence
of a BCMW-power for any periodic element.

For the periodic element $g$ in Definition~\ref{def:BCMW},
$q$ is the smallest positive integer such that $g^q$ belongs
to the cyclic group $\myangle{\Delta}$ up to conjugacy.
If $g$ is C-tight, then $g^q$ belongs to $\myangle{\Delta^m}$.
This notion is irrelevant to Theorem~\ref{thm:CMW}, but will
be used later for Theorem~\ref{thm:pa-cy}.

\begin{example}\label{eg:periodic}
Table~\ref{ta:ex} shows $\INF(\delta^k)$ and $\INF(\epsilon^k)$
for some braid index $n$
in $B_n^\Artinbr$, the $n$-braid group with the Artin Garside structure,
and in $B_n^\BKLbr$, the $n$-braid group with the BKL Garside structure.
In the table, $(*)$ indicates that the periodic braid is
P-minimal and $(**)$ indicates that it is P-minimal and C-tight.
We know that $\Delta^2$ and $\delta^n$ are the smallest
central power of $\Delta$ and $\delta$ in $B_n$.
Observe the following.

\begin{itemize}
\item
Let $G=B_n^\Artinbr$.
Since $\Delta$ is the Garside element
and $\delta^n = \epsilon^{n-1} = \Delta^2$,
$$
\INF(\delta)=2/n
\quad\mbox{and}\quad
\INF(\epsilon)=2/(n-1).
$$
Therefore, a $\delta$-type periodic $n$-braid is C-tight
if $n$ is odd, and
an $\epsilon$-type periodic $n$-braid is C-tight
if $n$ is even.

\item
Let $G=B_n^\BKLbr$.
Note that $\delta$ is the Garside element.
The super summit set of $\delta^k$ is $\{\delta^k\}$ for all integers $k$,
hence the conjugacy search problem is easy for $\delta$-type periodic braids.
Hence we are interested only in $\epsilon$-type periodic braids.
Since $\epsilon^{n-1}=\delta^n$,
$$
\INF(\epsilon)=n/(n-1).
$$
Since $n$ is relatively prime to $n-1$, every $\epsilon$-type periodic braid
is C-tight.
\end{itemize}
\end{example}

\begin{table}
$$\def\arraystretch{1.1}
\begin{array}{|c||l|l||l|l|}\hline
k
  & \textstyle\INF(\delta^k)~\mbox{in $B_{9}^\Artinbr$}
      \atop \textstyle\INF(\epsilon^k)~\mbox{in $B_{10}^\Artinbr$}
  & \textstyle\INF(\delta^k)~\mbox{in $B_{10}^\Artinbr$}
      \atop \textstyle\INF(\epsilon^k)~\mbox{in $B_{11}^\Artinbr$}
  & ~\INF(\epsilon^k)~\mbox{in $B_{10}^\BKLbr$}
  & ~\INF(\epsilon^k)~\mbox{in $B_{11}^\BKLbr$}\\\hline
1 & \qquad\frac29
  & (*)~{}~\frac{2}{10}=\frac15
  & (**)\,\frac{10}{9}=1+\frac19
  & (**)\,\frac{11}{10}=1+\frac{1}{10} \\\hline
2 & \qquad\frac 49
  & \qquad\frac{4}{10}=\frac25
  & \qquad\frac{20}{9}=2+\frac29
  & (**)\,\frac{22}{10}=\frac{11}5=2+\frac{1}{5}\\\hline
3 & \qquad\frac 69=\frac23
  & \qquad\frac{6}{10}=\frac35
  & (**)\,\frac{30}{9}=\frac{10}3=3+\frac13
  & \qquad\frac{33}{10}=3+\frac{3}{10} \\\hline
4 & \qquad\frac 89
  & \qquad\frac{8}{10}=\frac45
  & \qquad\frac{40}{9}=4+\frac49
  & \qquad\frac{44}{10}=\frac{22}5=4+\frac{2}{5} \\\hline
5 & (**)\,\frac {10}9=1+\frac19
  & (*)~{}~\frac{10}{10}=1
  & \qquad\frac{50}{9}=5+\frac59
  & (**)\,\frac{55}{10}=\frac{11}2=5+\frac{1}{2}\\\hline
6 & (**)\,\frac {12}9=\frac43=1+\frac13
  & (**)\,\frac{12}{10}=\frac65=1+\frac15
  & \qquad\frac{60}{9}=\frac{20}3=6+\frac23
  & \qquad\frac{66}{10}=\frac{33}5=6+\frac{3}{5} \\\hline
7 & \qquad\frac{14}9=1+\frac59
  & \qquad\frac{14}{10}=\frac75=1+\frac25
  & \qquad\frac{70}{9}=7+\frac79
  & \qquad\frac{77}{10}=7+\frac{7}{10} \\\hline
8 & \qquad\frac{16}9=1+\frac 79
  & \qquad\frac{16}{10}=\frac85=1+\frac35
  & \qquad\frac{80}{9}=8+\frac89
  & \qquad\frac{88}{10}=\frac{44}5=8+\frac{4}{5} \\\hline
9 & (**)\,\frac{18}{9}=2
  & \qquad\frac{18}{10}=\frac95=1+\frac45
  & (**)\,\frac{90}{9}=10
  & \qquad\frac{99}{10}=9+\frac{9}{10} \\\hline
\end{array}
$$
\caption{In the table, $(*)$ indicates that the periodic braid is
P-minimal and $(**)$ indicates that it is P-minimal and C-tight.}
\label{ta:ex}
\end{table}

As mentioned earlier, it follows from Theorem~\ref{thm:CMW} that
every periodic element has a BCMW-power.
We prove this in Theorem~\ref{thm:LL-CMW} in a different way
using only the properties of the invariants $\INF(\cdot)$,
$\SUP(\cdot)$ and $\LEN(\cdot)$, and the stable super summit set.
A new fact added by Theorem~\ref{thm:LL-CMW} is that
the exponent of a BCMW-power of a periodic element $g$
is completely determined only by $\INF(g)$ and $m$, further,
it can be computed easily from those values.
Before we go into the theorem, we show necessary lemmas.

\begin{lemma}\label{lem:per-equiv}
An element $g$ in $G$ is periodic if and only if\/
$\LEN(g)=0$, that is, $\INF(g)=\SUP(g)$.
\end{lemma}

\begin{proof}
Let $\INF(g)=p/q$, where $p$ and $q$ are relatively prime and $q\ge 1$.

Suppose that $g$ is periodic.
Since $g^k=\Delta^l$ for some integers $l$ and $k\ge 1$,
$$
\LEN(g)=\frac1k\cdot\LEN(g^k)=\frac1k\cdot\LEN(\Delta^l)=\frac1k\cdot 0=0.
$$

Conversely, assume that $\LEN(g)=0$, that is, $\INF(g)=\SUP(g)=p/q$.
Then
$$
\infs(g^q)=\lfloor \INF(g^q)\rfloor
=\lfloor q\cdot\INF(g)\rfloor =p
\quad\mbox{and}\quad
\sups(g^q)=\lceil \SUP(g^q)\rceil
=\lceil q\cdot\SUP(g)\rceil = p.
$$
Therefore $g^q$ is conjugate to $\Delta^p$, hence $g$ is periodic.
\end{proof}

\begin{lemma}\label{lem:per-elt}
Let $g$ be a periodic element of\/ $G$ with $\INF(g)=p/q$ for relatively prime
integers $p$ and $q$ with $q\ge 1$. Then the following hold.
\begin{itemize}
\item[(i)]
For all integers $n$,
$\INF(g^n) = n\cdot\INF(g)$.

\item[(ii)]
$\lens(g^k)$ is $0$ if\/ $k\equiv 0\bmod q$ and $1$ otherwise.
In particular, $[g^k]^S =[g^k]^U$ for all $k$.

\item[(iii)]
$q=1$ if and only if $g$ is conjugate to $\Delta^p$.
\end{itemize}
\end{lemma}

\begin{proof}
\smallskip\noindent
(i) \ \
We know that $\INF(g^n) = n\cdot\INF(g)$ holds if $n$ is nonnegative.
Let $n$ be negative, then $n=-k$ for some positive integer $k$.
Since $\INF(g)=\SUP(g)$ by Lemma~\ref{lem:per-equiv}
and $\INF(h^{-1})=-\SUP(h)$ for all $h\in G$,
we have
$$\INF(g^n) =\INF((g^{-1})^k)
=k\cdot\INF(g^{-1})
=k\cdot(-\SUP(g))=n\cdot\INF(g).
$$

\smallskip\noindent
(ii) \ \
For any positive integer $k$,
$$
\lens(g^{-k})
=\lens(g^k)
=\sups(g^k)-\infs(g^k)
=\lceil \SUP(g^k)\rceil -\lfloor\INF(g^k)\rfloor
=\lceil kp/q\rceil -\lfloor kp/q\rfloor.
$$
Therefore, for any integer $k$, $\lens(g^k)$ is either 0 or 1,
and $\lens(g^k)=0$ if and only if $kp/q$ is an integer.
Because $p$ and $q$ are relatively prime,
$kp/q$ is an integer if and only if $k\equiv 0\bmod q$.

\smallskip\noindent
(iii) \ \
Let $q=1$. Then $\INF(g)=\SUP(g)=p$ by Lemma~\ref{lem:per-equiv}.
Since $\infs(g)=\myfloor{\INF(g)}=p=\myceil{\SUP(g)}=\sups(g)$,
the element $g$ is conjugate to $\Delta^p$.
The converse is obvious.
\end{proof}

\begin{lemma}\label{lem:Primi_equiv}
Let $g$ be a periodic element of\/ $G$ with $\INF(g)=p/q$ for relatively
prime integers $p$ and $q\ge2$.
Then the following conditions are equivalent.
\begin{itemize}
\item[(i)]
$g$ is P-minimal (i.e. $p\equiv 1\bmod q$).

\item[(ii)]
Every element $h\in[g]^{St}$ is of the form $\Delta^u a$
for an integer $u$ and a simple element $a$ such that
the simple element $a$ is minimal in the set
$S = \{ \Delta^{-\inf(h^n)}h^n : n\in\Z\}\setminus\{e\}$
under the order relation $\le_R$.

\item[(iii)]
Every element $h\in [g]^{St}$ is of the form $\Delta^u a$
for an integer $u$ and a simple element $a$
such that $\tau^{(q-1)u}(a)\,\tau^{(q-2)u}(a)\cdots \tau^{u}(a)\,a = \Delta$.

\item[(iv)]
$g$ is conjugate to an element of the form $\Delta^u a$
for an integer $u$ and a simple element $a$
such that $\tau^{(q-1)u}(a)\,\tau^{(q-2)u}(a)\cdots \tau^{u}(a)\,a = \Delta$.
\end{itemize}
\end{lemma}

\begin{proof}
We prove the equivalences by showing the implications
(i) $\Rightarrow$ (iii) $\Rightarrow$  (iv) $\Rightarrow$ (i)
and (ii) $\Leftrightarrow$  (iii).
Note that the implication (iii) $\Rightarrow$ (iv) is obvious.

\smallskip\noindent
(i) $\Rightarrow$ (iii)\ \
Suppose that $\INF(g)=(uq+1)/q=u+1/q$ for an integer $u$.
Because $g$ is periodic, $\SUP(g)=\INF(g)=u+1/q$ by Lemma~\ref{lem:per-equiv}.
Let $h$ be an element of the stable super summit set of $g$.
Then
$$
\begin{array}{ll}
\inf(h)=\infs(g)=\lfloor \INF(g)\rfloor=u,\quad
&\inf(h^q)=\infs(g^q)=\lfloor q\cdot \INF(g)\rfloor=uq+1,\\
\sup(h)=\sups(g)=\lceil \SUP(g)\rceil=u+1,\quad
&\sup(h^q)=\sups(g^q)=\lceil q\cdot \SUP(g)\rceil=uq+1.
\end{array}
$$
Therefore, $h=\Delta^u a$ for $a\in\D\backslash\{ e, \Delta\}$
and $h^q=\Delta^{uq+1}$.
Note that
$$
\Delta^{uq+1}=h^q=(\Delta^u a)\cdots(\Delta^u a)
=\Delta^{uq}\, \tau^{(q-1)u}(a)\,\tau^{(q-2)u}(a)\cdots \tau^{u}(a)\,a.
$$
This  implies that
$\Delta=\tau^{(q-1)u}(a)\,\tau^{(q-2)u}(a)\cdots \tau^{u}(a)\,a$.

\smallskip\noindent
(iv) $\Rightarrow$ (i)\ \
Suppose that $g$ is conjugate to $h=\Delta^u a$ for $a\in\D\backslash\{ e, \Delta\}$
satisfying
$$\tau^{(q-1)u}(a)\,\tau^{(q-2)u}(a)\cdots \tau^{u}(a)\,a = \Delta.$$
Then $h^q=\Delta^{uq+1}$,
hence one has $\INF(h)=(1/q)\cdot\INF(h^q)=(uq+1)/q$.
Since $p/q = \INF(g)=\INF(h)=(uq+1)/q$, one has $p=uq+1\equiv 1\bmod q$.

\smallskip\noindent
(ii) $\Rightarrow$ (iii)\ \
For $j\ge 1$, let $a_j$ be the positive element defined by
$$
a_j=\tau^{(j-1)u}(a)\,\tau^{(j-2)u}(a)\cdots \tau^{u}(a)\,a.
$$
Note that $h^j=(\Delta^u a)^j=\Delta^{ju} a_j$, hence
$\len(a_j)=\len(h^j)=\lens(g^j)$ and it is 0 if $j\equiv 0\bmod q$
and 1 otherwise.
In particular,
$\len(a_q)=\sup(a_q)-\inf(a_q)=0$, hence
$$
\inf(a_q)=\sup(a_q)\ge \sup(a)\ge 1.
$$
Since $\inf(a_1)=\inf(a)=0$ and
$0\le \inf(a_{j+1})\le \inf(a_j)+1$ for all $j\ge 1$,
there exists $2\le j\le q$ such that $\inf(a_j)=1$.
Let $k$ be the smallest positive integer such that $\inf(a_k)=1$,
then $2\le k\le q$.
Since $\inf(a_{k-1})=0$ and $\len(a_{k-1})=1$,
$a_{k-1}$ is a simple element.
Since
$$
\Delta\le_L a_k
=\tau^{(k-1)u}(a)\cdots \tau^{u}(a)\,a
=\tau^u\bigl(\tau^{(k-2)u}(a)\cdots \tau^{u}(a)a\bigr)\,a
=\tau^u(a_{k-1}) a,
$$
we have $a=a'a''$ for some simple elements $a'$ and $a''$ such that
$\tau^u(a_{k-1})a'=\Delta$, that is, $a_k=\Delta a''$.
Since $\inf(a_{k-1})=0$, we have $a'\ne e$, hence
$a''\ne a$.
Note that $a''\le_R a$ and $h^k=\Delta^{uk+1}a''$.
If $a''\ne e$, then it contradicts the minimality of $a$ under $\le_R$.
Therefore $a''=e$, hence $a_k=\Delta$.
Notice that $\len(a_j)=0$ if and only if $j\equiv 0\bmod q$.
Since $2\le k\le q$ by the construction,
one has $k=q$ and hence
$$
\tau^{(q-1)u}(a)\,\tau^{(q-2)u}(a)\cdots \tau^{u}(a)\,a=a_q=a_k=\Delta
$$
as desired.

\smallskip\noindent
(iii) $\Rightarrow$ (ii)\ \
It is obvious since $S=\{e\}\cup\{ \tau^{(k-1)u}(a)\,\tau^{(k-2)u}(a)\cdots
\tau^{u}(a)\,a: k=1,\ldots,q-1\}$.
\end{proof}

\begin{lemma}\label{lem:SubGp_equiv}
Let $g$ be a periodic element of\/ $G$ with $\INF(g)=p/q$ for relatively
prime integers $p$ and $q\ge 1$, and let $r$ be a nonzero integer.
Then, the elements $\bar g$ and $\bar g^r$ generate the same cyclic subgroup of $G_\Delta$
if and only if $r$ is relatively prime to both
$q$ and $m/\gcd(p,m)$.
\end{lemma}

\begin{proof}
Let $\langle\bar g\rangle$ be the cyclic subgroup
of $G_\Delta$ generated by $\bar g$.
For positive integers $k$,
let $\Z/k\Z$ denote the additive group of residue classes modulo $k$.
Define a map $\phi:\langle\bar g\rangle \to \Z/(qm)\Z$ by
$$
\phi(\bar g^n)=q\cdot\INF(g^n)=q\cdot (np/q)=pn \bmod {qm},
\qquad n\in\Z.
$$
We will first show that $\phi$ is an injective homomorphism.
$\phi$ is well-defined because $q\cdot\INF(\Delta^m)=qm\equiv 0\bmod qm$.
$\phi$ is a homomorphism because $\INF(g^n)=n\cdot\INF(g)$ for all integers $n$.
If $\phi(\bar g^n)\equiv 0\bmod {qm}$,
then $\INF(g^n)=km$ for some integer $k$.
Because $g^n$ is periodic, we have $\SUP(g^n)=\INF(g^n)=km$,
hence $\infs(g^n)=\sups(g^n)=km$.
Therefore $g^n$ is conjugate to $\Delta^{km}$.
Because $\Delta^m$ is central, one has $g^n=\Delta^{km}$.
This shows that $\phi$ is injective.

Let $r$ be a nonzero integer.
Since $\phi(\bar g^n)=pn \bmod{qm}$ for all integers $n$,
the images of $\langle\bar g\rangle$ and $\langle\bar g^r\rangle$ under $\phi$
are generated by the residue classes of $\gcd(p,qm)$ and $\gcd(pr,qm)$,
respectively.
Therefore it suffices to show that
$\gcd(p,qm)=\gcd(pr,qm)$ if and only if
$r$ is relatively prime to both $q$ and $m/\gcd(p,m)$.

For nonzero integers $a$, $b$ and $c$, the following equality holds:
$$
\gcd(ab,c)=\gcd(a,c)\cdot\gcd\left(b,{c}/{\gcd(a,c)}\right).
$$
Note that $\gcd(p,qm)=\gcd(p,m)$ because $p$ and $q$ are relatively prime.
Applying the above equality to $a=p$, $b=r$ and $c=qm$,
\begin{eqnarray*}
\gcd(pr,qm)
&=&\gcd(p,qm)\cdot\gcd (r,qm/\gcd(p,qm) )\\
&=&\gcd(p,qm)\cdot\gcd (r,qm/\gcd(p,m) ).
\end{eqnarray*}
Therefore, $\gcd(pr,qm)=\gcd(p,qm)$ if and only if
$\gcd(r, qm/\gcd(p,m))=1$.
Since $qm/\gcd(p,m)=q\cdot (m/\gcd(p,m))$ and $m/\gcd(p,m)$ is
an integer, $\gcd(r, qm/\gcd(p,m))=1$
if and only if
$r$ is relatively prime to both $q$ and $m/\gcd(p,m)$.
\end{proof}

\begin{theorem}\label{thm:LL-CMW}
Let $g$ be a periodic element of a Garside group\/ $G$ with $\INF(g)=p/q$ for relatively
prime integers $p$ and $q$ with $q\ge 1$.
A power $g^r$ is a BCMW-power of $g$ if and only if\/
$pr\equiv 1 \bmod q$ and $r$ is relatively prime to $m/\gcd(p,m)$.
In particular, if\/ $q=1$ then $g$ itself is a BCMW-power, and
if\/ $q\ge 2$ then there is an integer $r$ with $1\le r<qm$ such that
$g^r$ is a BCMW-power of $g$.
\end{theorem}

\begin{proof}
By Lemma~\ref{lem:per-elt}, one has $\INF(g^r)=pr/q$ for any integer $r$.

Suppose that $g^r$ is a BCMW-power of $g$.
Since $\myangle{\bar g}=\myangle{\bar g^r}$ in $G_\Delta$,
$r$ is relatively prime to $q$ and $m/\gcd(p,m)$
by Lemma~\ref{lem:SubGp_equiv}.
Since $q$ is relatively prime to $p$ by assumption,
$q$ is relatively prime to $pr$.
Therefore, $pr\equiv 1 \bmod q$ because $g^r$ is P-minimal.

Conversely, suppose that
$pr\equiv 1 \bmod q$ and $r$ is relatively prime to $m/\gcd(p,m)$.
Since $pr$ and $q$ are relatively prime, so are $r$ and $q$.
Hence $\myangle{\bar g}=\myangle{\bar g^r}$ in $G_\Delta$ by Lemma~\ref{lem:SubGp_equiv}.
Since $pr\equiv 1 \bmod q$, $g^r$ is P-minimal.
Therefore $g^r$ is a BCMW-power of $g$.

\smallskip

If $q=1$, then $g$ itself is obviously a BCMW-power
because it is P-minimal.

\smallskip

We now show the existence of a BCMW power $g^r$ with $1\le r<qm$
for $q\ge 2$.
Since $p$ and $q$ are relatively prime, there exists an integer $r_0$ such that
$1\le r_0 < q$ and
$$
pr_0\equiv 1\bmod q.
$$
In particular, $r_0$ is relatively prime to $q$.

Suppose that all of the prime divisors of $m/\gcd(p,m)$
are also divisors of $q$.
In this case, let $r=r_0$.
Then $r$ is relatively prime to $m/\gcd(p,m)$
since it is relatively prime to $q$.
Moreover $1\le r < q\le qm$ and $pr\equiv 1\bmod q$, as desired.

Let $p_1,\ldots,p_k$ be all of the distinct prime divisors
of $m/\gcd(p,m)$ which do not divide $q$.
For each $p_i$, take an integer $r_i$ relatively prime to $p_i$.
By the Chinese remainder theorem, there exists an integer $r$
with $1\le r<q p_1 p_2 \cdots p_k$ such that
$$
r\equiv r_0\bmod {q},\qquad
r\equiv r_i\bmod {p_i}\quad\mbox{for $i=1,\ldots,k$}.
$$
By the construction, $r$ is relatively prime to $m/\gcd(p,m)$
and $pr\equiv pr_0\equiv 1\bmod q$.
Moreover $1\le r<q p_1 p_2 \cdots p_k
\le q\cdot m/\gcd(p,m)\le qm$, as desired.
\end{proof}

Combining with Lemmas~\ref{lem:Primi_equiv} and \ref{lem:SubGp_equiv},
the above result implies the following.

\begin{corollary}\label{cor:LL-CMW}
Let $g$ be a periodic element of\/ $G$ with $\INF(g)=p/q$ for relatively
prime integers $p$ and $q$ with $q\ge 2$.
Then there exists an integer $r$ such that the following hold.
\begin{itemize}
\item
$g^r$ generates the same cyclic subgroup of $G_\Delta$ as $g$ does.

\item
$g^r$ is conjugate to $\Delta^u a$ for a simple element
$a\in\D\backslash\{e, \Delta\}$ such that
$$
\tau^{(q-1)u}(a)\,\tau^{(q-2)u}(a)\cdots \tau^{u}(a)\,a=\Delta.
$$
\end{itemize}
The exponent $r$ depends only on the integers $p$, $q$ and $m$,
and it is characterized by the property that
$r$ is relatively prime to $m/\gcd(p,m)$ and $pr\equiv 1\bmod q$.
In particular, one can take $r$ such that $1\le r<qm$.
\end{corollary}

\begin{remark}
Theorem~\ref{thm:CMW} shows the structure of finite subgroups of $G_\Delta$.
In particular, every finite subgroup of $G_\Delta$ is abelian.
Once we know that finite subgroups of $G_\Delta$ are abelian,
then it is easy to prove Theorem~\ref{thm:CMW}
using only the properties of the invariants $\INF(\cdot)$, $\SUP(\cdot)$ and
$\LEN(\cdot)$, and the stable super summit set.
However, it looks difficult to prove that finite subgroups
of $G_\Delta$ are abelian without considering the action of $G$
on Bestvina's normal form complex as done in~\cite{Bes99} and in~\cite{CMW04}.
\end{remark}

For $g, h, x\in G$ and $n\in\Z$, $x^{-1}gx=h$ implies $x^{-1}g^nx=h^n$.
But the converse direction is not true in general.
Even though $g^n$ is conjugate to $h^n$,
$g$ is not necessarily conjugate to $h$.
Interestingly, the converse direction holds for BCMW-powers.

\begin{proposition}\label{prop:BCMW_power}
Let $g$ be a periodic element of\/ $G$ and $g^r$ a BCMW-power.
Then for elements $h$ and $x$ of\/ $G$,
$x$ conjugates $g$ to $h$
if and only if it conjugates $g^r$ to $h^r$.
Therefore, the conjugacy decision problem and
the conjugacy search problem for $(g,h)$
are equivalent to those for $(g^r,h^r)$,
and the centralizer of\/ $g$ in $G$ is the same as
the centralizer of\/ $g^r$ in $G$.
\end{proposition}

\begin{proof}
If $x^{-1}gx=h$, then it is obvious that $x^{-1}g^r x=h^r$.
Conversely, suppose that $x^{-1}g^r x=h^r$.
We claim that there exist integers $s$ and $t$ such that
$$
g^{rs}=\Delta^{mt}g\quad\mbox{and}\quad
h^{rs}=\Delta^{mt}h.
$$

Since $\langle \bar g\rangle=\langle\bar g^r\rangle$ in $G_\Delta$,
there exist integers $s$ and $t$ such that
$g^{rs}=\Delta^{mt}g$.
Since $h^r$ and $g^r$ are conjugate, $h^r$ is periodic
(hence $h$ is periodic) and $\INF(h^r)=\INF(g^r)$.
Since $g^{rs-1}=\Delta^{mt}$ and $h$ is periodic,
$$
\SUP(h^{rs-1})=\INF(h^{rs-1})
=\frac{rs-1}r\cdot \INF(h^r)
=\frac{rs-1}r\cdot \INF(g^r)
=\INF(g^{rs-1})=mt.
$$
Since $mt$ is an integer, we have $\infs(h^{rs-1})=\sups(h^{rs-1})=mt$,
hence $h^{rs-1}$ is conjugate to $\Delta^{mt}$.
Since $\Delta^{mt}$ is central, $h^{rs-1}=\Delta^{mt}$.
Therefore $h^{rs}=\Delta^{mt}h$.

Since $x^{-1}g^r x=h^r$, we have
$x^{-1}(\Delta^{mt}g)x=x^{-1}g^{rs}x=h^{rs}=\Delta^{mt}h$.
Since $\Delta^{mt}$ is central,
it follows that $x^{-1}g x=h$.
\end{proof}

\subsection{Super summit sets of P-minimal and C-tight periodic elements}
Recall the definition of partial cycling introduced by Birman, Gebhardt
and Gonzalez-Meneses~\cite{BGG06b}.

\begin{definition}
Let $g=\Delta^u a_1a_2\cdots a_\ell \in G$ be in the normal form.
Let $b\in\D$ be a prefix of $a_1$, i.e. $a_1=ba_1'$ for a simple element $a'_1$.
The conjugation
$$
\tau^{-u}(b)^{-1} g \tau^{-u}(b) = \Delta^u a_1'a_2\cdots a_\ell\tau^{-u}(b)
$$
is called a \emph{partial cycling} of $g$ by $b$.
\end{definition}

\begin{remark}
For a ultra summit element $g=\Delta^u a_1a_2\cdots a_\ell \in G$ in the normal form
with $\ell > 0$,
a simple element $b\neq e$ is called a {\em minimal simple element} for $g$
with respect to $[g]^U$ if $b^{-1}gb\in [g]^U$ and no proper prefix of $b$ satisfies
this property.
When partial cycling was defined at first in~\cite{BGG06b},
it was used only for conjugating $g$
by a minimal simple element, a special prefix of $\tau^{-u}(a_1)$.
Unlike the previous usage, we deal with partial cycling which conjugates $g$
by any prefix of $\tau^{-u}(a_1)$.
\end{remark}

We now establish the main result (Theorem~\ref{thm:pa-cy}) of this section
that the super summit set of a P-minimal, C-tight
periodic element is closed under any partial cycling.
Joining Theorem~\ref{thm:LL-CMW} that every periodic element has a BCMW-power,
Theorem~\ref{thm:pa-cy} yields that the super summit set of a BCMW-power of a C-tight
periodic element is closed under any partial cycling.
Notice that, without one of the conditions in Theorem~\ref{thm:pa-cy},
the super summit set is not necessarily closed under partial cycling.
We will show this in Example~\ref{eg:conditions}.

\begin{theorem}\label{thm:pa-cy}
Let $g$ be a periodic element of a Garside group $G$.
If $g$ is P-minimal and C-tight, then
the following conditions are equivalent for an element $h$ conjugate to $g$.
\begin{itemize}
\item[(i)] $\inf(h)=\infs(g)$.
\item[(ii)] $h\in [g]^{S}$, that is, $\inf(h)=\infs(g)$ and $\sup(h)=\sups(g)$.
\item[(iii)] $h\in[g]^{St}$, that is, $\inf(h^k)=\infs(g^k)$
and $\sup(h^k)=\sups(g^k)$ for all $k\ge 1$.
\end{itemize}
In particular, $[g]^{S}$ is closed under any partial cycling.
\end{theorem}

\begin{proof}
Let $\INF(g)=p/q$ for relatively prime integers $p$ and $q$ with $q\ge 1$.
The claim is trivial when $q=1$, because $g$ is conjugate to $\Delta^p$.
Thus we may assume that $q\ge 2$.

Because the implications (iii) $\Rightarrow$ (ii) $\Rightarrow$ (i)
are obvious, we will show that (i) $\Rightarrow$ (iii).
Suppose that $h$ is conjugate to $g$ and $\inf(h)=\infs(g)$.

Because $g$ is P-minimal and C-tight, $p=uq+1=ml$ for some integers $u$ and $l$.
Therefore
\begin{eqnarray*}
\INF(h)&=&\SUP(h)=p/q=u+1/q,\\
\infs(h^k)&=& \lfloor k\INF(h)\rfloor = \lfloor k(uq+1)/q\rfloor=ku+\lfloor k/q\rfloor,\\
\sups(h^k)&=& \lceil k\SUP(h)\rceil = \lceil k(uq+1)/q\rceil = ku+\lceil k/q\rceil
\end{eqnarray*}
for all integers $k\ge 1$.
Because $\infs(h^q)=\sups(h^q)=uq+1=ml$,
$h^q$ is conjugate to $\Delta^{ml}$.
Because $\Delta^{ml}$ is central,
$$
h^q=\Delta^{ml}=\Delta^p=\Delta^{uq+1}.
$$
Because $\inf(h)=\infs(h)=u$, there exists a positive element $a$ such that
$$
h=\Delta^u a.
$$
Because $\sup(h)\ge \sups(h)=u+1$, $a$ is not the identity.
Let $\psi$ denote $\tau^{u}$. Then
$$
h^k=\Delta^{ku}\psi^{k-1}(a)\psi^{k-2}(a)\cdots\psi(a) a
$$
for all integers $k\ge 2$.
Because $h^q=\Delta^{qu+1}$,
$$
\Delta^{qu}\cdot\Delta=\Delta^{qu+1}=h^q
=\Delta^{qu}\psi^{q-1}(a)\psi^{q-2}(a)\cdots\psi(a) a.
$$
Therefore $\Delta=\psi^{q-1}(a)\psi^{q-2}(a)\cdots\psi(a) a$.
In particular,
$\psi^{k-1}(a)\psi^{k-2}(a)\cdots a \in \D\backslash\{ e, \Delta\}$
for all integers $k$ with $1\le k < q$.
Therefore
$$
\inf(h^k)=uk=\infs(h^k)
\quad\mbox{and}\quad
\sup(h^k)=uk+1=\sups(h^k)
\qquad\mbox{for $k=1,\ldots,q-1$}.
$$
Because $h^q=\Delta^{ml}$, this proves that
$h$ belongs to the stable super summit set.

\medskip
Now, let us show that $[g]^S$ is closed under any partial cycling.
Let $h'$ be the result of an arbitrary partial cycling
of an element $h$ in $[g]^S$.
Partial cycling does not decrease the infimum by definition.
Therefore $\inf(h')=\inf(h)=\infs(g)$, whence
$h'$ belongs to $[g]^S$.
\end{proof}

For $g\in G$, Garside~\cite{Gar69} called the set
$\{h\in[g]: \inf(h)=\infs(h)\}$ the \emph{summit set} of $g$.
Theorem~\ref{thm:pa-cy} shows that for P-minimal, C-tight
periodic elements, the notions of summit set, super summit set
and stable super summit set are all equivalent.
We already know that, for a periodic element, its summit length is at most 1, hence
its ultra summit set is nothing more than its super summit set.

The following example shows that the conditions in
Theorem~\ref{thm:pa-cy} that the periodic element $g$
is P-minimal and C-tight are necessary for the conclusion.

\begin{example}\label{eg:conditions}
Recall that $B_n^\Artinbr$ and $B_n^\BKLbr$ denote
the $n$-braid group with the Artin Garside structure
and the BKL Garside structure, respectively.
We will observe the following.
\begin{itemize}
\item $\epsilon_{(5)}=\sigma_4\sigma_3\sigma_2\sigma_1\sigma_1\in B_5^\Artinbr$
is P-minimal but not C-tight because $\INF(\epsilon_{(5)})=1/2$.
\item $\epsilon_{(6)}=\sigma_5\sigma_4\sigma_3\sigma_2\sigma_1\sigma_1\in B_6^\Artinbr$
is C-tight but not P-minimal because $\INF(\epsilon_{(6)})=2/5$.

\item
$\epsilon_{(6)}^3=\delta_{(6)}^3[4,3,2,1]\in B_6^\BKLbr$
is C-tight but not P-minimal
because $\INF(\epsilon_{(6)}^3)=18/5$.

\item
For all of these examples, their super summit sets
are not closed under partial cycling
and different from their stable super summit sets.
\end{itemize}

Note that $\epsilon_{(n)}^{n-1}=\Delta^2=\delta_{(n)}^n$,
hence $\INF(\epsilon_{(n)})=2/(n-1)$ in $B_n^\Artinbr$
and $\INF(\epsilon_{(n)})=n/(n-1)$ in $B_n^\BKLbr$.
Therefore
$$
\begin{array}{ll}
\INF(\epsilon_{(5)}) = \frac24=\frac12
& \mbox{in $B_5^\Artinbr$},\\
\INF(\epsilon_{(6)}) = \frac25
& \mbox{in $B_6^\Artinbr$},\\
\INF(\epsilon_{(6)}^3) = 3\cdot\frac{6}{5}=\frac{18}5=3+\frac35
& \mbox{in $B_6^\BKLbr$}.
\end{array}
$$

Consider the elements
$g_1=\sigma_1\sigma_4\sigma_3\sigma_2\sigma_1$ in $[\epsilon_{(5)}]^S$
and $g_2=\sigma_1\sigma_5\sigma_4\sigma_3\sigma_2\sigma_1$
in $[\epsilon_{(6)}]^S$, under the Artin Garside structure.
The partial cycling on $g_1$ and $g_2$ by $\sigma_1$ yields
$\epsilon_{(5)}=(\sigma_4\sigma_3\sigma_2\sigma_1)\sigma_1$ and
$\epsilon_{(6)}=(\sigma_5\sigma_4\sigma_3\sigma_2\sigma_1)\sigma_1$, respectively.
Neither $\epsilon_{(5)}$ nor $\epsilon_{(6)}$ is a super summit element,
because they have canonical length 2.
Hence, $[\epsilon_{(5)}]^S$ and $[\epsilon_{(6)}]^S$
are not closed under partial cycling.

The normal forms of $g_1^2$ and $g_2^2$ are as in the right hand sides
in the following equations:
\begin{eqnarray*}
g_1^2
&=&(\sigma_1\sigma_4\sigma_3\sigma_2\sigma_1)
   (\sigma_1\sigma_4\sigma_3\sigma_2\sigma_1)
  =(\sigma_1\sigma_4\sigma_3\sigma_2\sigma_1\sigma_4\sigma_3\sigma_2)
   \cdot(\sigma_1\sigma_2),\\
g_2^2
&=&(\sigma_1\sigma_5\sigma_4\sigma_3\sigma_2\sigma_1)
   (\sigma_1\sigma_5\sigma_4\sigma_3\sigma_2\sigma_1)
  =(\sigma_1\sigma_5\sigma_4\sigma_3\sigma_2\sigma_1\sigma_5\sigma_4\sigma_3\sigma_2)
   \cdot(\sigma_1\sigma_2).
\end{eqnarray*}
In particular, $g_1^2$ and $g_2^2$ have canonical length 2, hence
they do not belong to their super summit sets.
Hence $[\epsilon_{(5)}]^S\ne[\epsilon_{(5)}]^{St}$
and $[\epsilon_{(6)}]^S\ne[\epsilon_{(6)}]^{St}$.

\smallskip
Now we consider $\epsilon_{(6)}^3 \in B_6^\BKLbr$.
Clearly, $\epsilon_{(6)}^3$ belongs to its super summit set
because it has canonical length 1.
Note that
$$
\epsilon_{(6)}^3=\delta_{(6)}^3[4,3,2,1]=\delta_{(6)}^3[4,1][4,3,2]
=\delta_{(6)}^3[4,2][4,3][2,1].
$$
Partial cycling of $\epsilon_{(6)}^3$ by $[4,1]$ gives
$$
\delta_{(6)}^3 [4,3,2] \tau^{-3}([4,1])
= \delta_{(6)}^3 [4,3,2]\,[4,1].
$$
Note that $[4,3,2][4,1]$ is not a simple element,
hence $[\epsilon_{(6)}^3]^S$ is not closed under partial cycling.
Let $g_3$ be the result of partial cycling of $\epsilon_{(6)}^3$ by $[4,2]$,
that is,
$$
g_3
=\delta_{(6)}^3 [4,3][2,1] \tau^{-3}([4,2])
=\delta_{(6)}^3 [4,3][2,1] [5,1]
=\delta_{(6)}^3 [4,3][5,2,1].
$$
Then $g_3\in [\epsilon_{(6)}^3]^S$ because $[4,3][5,2,1]$ is a simple element.
On the other hand, the normal form of $g_3^2$ is as in the right hand side
of the following equation.
$$
g_3^2
= \delta_{(6)}^3 [4,3][5,2,1] \delta_{(6)}^3 [4,3][5,2,1]
= \delta_{(6)}^6 [6,1][5,4,2] [4,3][5,2,1]
= \delta_{(6)}^6 [6,1][5,4,3,2] \cdot [5,2,1]
$$
Because $g_3^2$ has canonical length 2,
it does not belong to its super summit set.
Hence $[\epsilon_{(6)}^3]^S\ne[\epsilon_{(6)}^3]^{St}$.
\end{example}

From Theorems~\ref{thm:LL-CMW} and~\ref{thm:pa-cy}, we have the following result.

\begin{corollary}\label{cor:pa-cy}
Let $g$ be a periodic element of a Garside group $G$.
If $g$ is C-tight, then the super summit set of a BCMW-power of $g$
is closed under any partial cycling.
\end{corollary}

\begin{proof}
Let $\INF(g)=p/q$ for relatively prime integers $p$ and $q$ with $q\ge 1$.
The claim is trivial if $q=1$, so we may assume that $q\ge 2$.

Let $g^r$ be an arbitrary BCMW-power of $g$.
Then $\INF(g^r)=pr/q$ (by Lemma~\ref{lem:per-elt}), where
$pr$ is relatively prime to $q$ (by Theorem~\ref{thm:LL-CMW}).
Since $g$ is C-tight, one has $p\equiv 0\bmod m$
and hence $pr\equiv 0\bmod m$, which means that $g^r$
is a C-tight periodic element.
On the other hand, $g^r$ is P-minimal (by definition),
hence $[g^r]^{S}$ is closed under partial cycling
by Theorem~\ref{thm:pa-cy}.
\end{proof}

\section{Super summit sets of $\epsilon$-type periodic braids}
\label{sec:SSS_of_e}

In this section we consider the BKL Garside structure for the braid group $B_n$.
In this Garside group, the braid $\delta\ (= a_{n,n-1}\cdots a_{3,2}a_{2,1})$
is the Garside element.
This means that the super summit set
of $\delta^k$ consists of a single element $\delta^k$ for all integers $k$.
Therefore, the conjugacy search problem for $\delta$-type periodic braids is easy.
But it does not hold for $\epsilon$-type periodic braids.
Notice that, for all integers $k$, $\len(\epsilon^k)\le 1$ and
$\epsilon^k\in [\epsilon^k]^S = [\epsilon^k]^U$ in $B_n^{[\BKL]}$.

The main results of this section are Proposition~\ref{prop:size}
and Proposition~\ref{prop:main}.
\begin{itemize}
\item
In Proposition~\ref{prop:size} we show that the size of the ultra summit set
of $\delta^u[k,k-1,\ldots,1]$ is at least $\C_k$, the $k$th Catalan number
if $u$ and $k$ satisfy some constraints.
Asymptotically, the Catalan numbers grow as
$$
\C_k=\frac1{k+1}{2k\choose k}\sim\frac{4^k}{k^{\frac32}\sqrt\pi}.
$$
This implies that the ultra summit set of $\epsilon^k$ in $B_n^\BKLbr$
is exponentially large with respect to $n$ for some $k$.

\item
In Proposition~\ref{prop:main} we show that
by applying polynomially many partial cyclings to an arbitrary super summit
element of $\epsilon^d$ for proper divisors $d$ of $n-1$, we obtain $\epsilon^d$.
\end{itemize}

\begin{proposition}\label{prop:size}
Let $\alpha$ be an $n$-braid such that
$$
\alpha=\delta^u[k,k-1,\ldots,1],\qquad
2\le k\le u\le\frac n2\quad\mbox{or}\quad 2\le k=u+1\le\frac{n}2.
$$
Then the cardinality of\/ $[\alpha]^U$ is at least
$\C_k$, the $k$th Catalan number.
\end{proposition}

\begin{proof}
From~\cite{BKL98}, there are $\C_k$ left divisors of $[k,k-1,\ldots,1]$.
Let $a_i$, $i=1,\ldots,\C_k$, denote them.
For each $i$, let $b_i$ be the simple element satisfying
$a_ib_i=[k,k-1,\ldots,1]$, hence
$$
\alpha=\delta^u[k,k-1,\ldots,1]=\delta^u a_ib_i,\qquad i=1,2,\ldots,\C_k.
$$
For $i=1,2,\ldots,\C_k$, let $\alpha_i$ be an element such that
$$
\alpha_i=b_i\alpha b_i^{-1}=b_i\delta^u a_i=\delta^u\tau^u(b_i)a_i.
$$

First, we show that each $\alpha_i$ is a ultra summit element.
Observe that, for all $i$,
$$
\tau^u(b_i)a_i
\le_L \tau^u(b_i)a_ib_i
\le_R \tau^u(a_ib_i) a_ib_i
= [u+k,\ldots,u+2,u+1]\,[k,\ldots,2,1].
$$
Notice that $2\le k\le u+1< u+k\le n$.
If $k=u+1$, then
$$
[u+k,\ldots,u+1]\,[k,\ldots,2,1]
=[u+k,\ldots,u+1]\,[u+1,\ldots,2,1]
=[u+k,\ldots,2,1].
$$
If $k\le u$, then the descending cycles $[u+k,\ldots,u+2,u+1]$ and
$[k,\ldots,2,1]$ are parallel.
In either case, $[u+k,\ldots,u+2,u+1]\,[k,\ldots,2,1]$ is a simple element.
Therefore $\tau^u(b_i)a_i$ is a simple element, hence
$\alpha_i$ belongs to the ultra summit set $[\alpha]^U$ for all $i$.

\smallskip
Next, we show that $\alpha_i$'s are all distinct.
Suppose that $\alpha_i=\alpha_j$ for some $i$ and $j$.
Because $\tau^u(b_i)a_i=\tau^u(b_j)a_j$,
$$
\tau^u(b_i^{-1}b_j)=a_ia_j^{-1}.
$$
Note that $a_ia_j^{-1}$ belongs to the subgroup
$\langle \sigma_1,\ldots,\sigma_{k-1}\rangle$
whereas $\tau^u(b_i^{-1}b_j)$ belongs to the subgroup
$\langle \sigma_{u+1},\ldots,\sigma_{u+k-1}\rangle$.
Because
$\langle \sigma_1,\ldots,\sigma_{k-1}\rangle
\cap \langle \sigma_{u+1},\ldots,\sigma_{u+k-1}\rangle
=\{e\}$, we obtain
$\tau^u(b_i^{-1}b_j)=a_ia_j^{-1}=e$.
Therefore $a_i=a_j$, hence $i=j$.
\end{proof}

\begin{figure}
\begin{tabular}{ccc}
\includegraphics{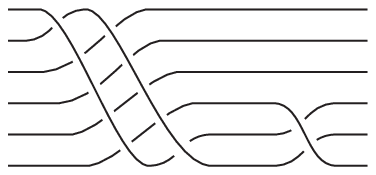} &\qquad&
\includegraphics{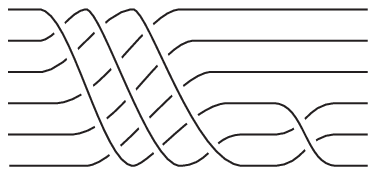} \\
(a)\ \ $\delta^2[3,2,1]\in B_6^\BKLbr$ &&
(b)\ \ $\delta^3[3,2,1]\in B_6^\BKLbr$
\end{tabular}
\caption{}
\label{fig:ex}
\end{figure}

\begin{remark}\label{rmk:size}
Let $\alpha\in B_{n}^{[\BKL]}$ be as in the above proposition.
Notice that if $k$ is proportional to $n$ such as $k\sim n/2$,
then the size of the ultra summit set of $\alpha$
is exponential in the braid index $n$.
Hence Proposition~\ref{prop:size} shows that the ultra summit sets of
the following braids are huge in $B_{n}^{[\BKL]}$.
\begin{itemize}
\item
If $k=u+1$, then $\alpha=\delta^{k-1}[k,k-1,\ldots,1]=\epsilon^{k-1}$.
(See Figure~\ref{fig:ex}~(a).)

\item If $k\le u$ and $u$ is a divisor of $n$,
then $\alpha$ is reducible.
(See Figure~\ref{fig:ex}~(b).)

\item If $k\le u$ and $u$ is not a divisor of $n$, then
$\alpha$ seems to be a pseudo-Anosov braid.
However, it would be beyond the scope of this paper to prove it.
\end{itemize}
\end{remark}

\begin{lemma}\label{lem:e^d}
For any positive divisor $d$ of\/ $n-1$,
the element $\epsilon^d$ is P-minimal and C-tight.
In particular, $[\epsilon^d]^{S}$ is closed under any partial cycling.
\end{lemma}

\begin{proof}
Notice that $n$ is the smallest positive integer such that $\delta^n$
is central.
Let $q=(n-1)/d$, hence $n=dq+1$.
Then $q$ and $n$ are relatively prime and
$$
\INF(\epsilon^d) = \frac{nd}{n-1}=\frac{n}{q}=\frac{dq+1}{q}.
$$
Therefore $\epsilon^d$ is P-minimal because $dq+1\equiv 1\bmod q$,
and C-tight because $dq+1=n\equiv 0\bmod n$.
By Theorem~\ref{thm:pa-cy}, the super summit set
$[\epsilon^d]^{S}$ is closed under partial cycling.
\end{proof}

The following lemma shows how to find a C-tight BCMW power of an arbitrary
$\epsilon$-type periodic braid.

\begin{lemma}\label{lem:exp_red}
Let $\alpha\in B_n$ be conjugate to $\epsilon^k$ for $k\ne 0$.
Let $d=\gcd(k,n-1)$, and $r$ and $s$ be integers such that $kr+(n-1)s=d$.
\begin{itemize}
\item[(i)] $\alpha^r$ is a C-tight BCMW-power of\/ $\alpha$,
being conjugate to $\delta^{-ns}\epsilon^d$.
\item[(ii)] Let $\alpha_1=\delta^{ns}\alpha^r$.
An $n$-braid $\gamma$ conjugates $\alpha$ to $\epsilon^k$
if and only if\/ it conjugates $\alpha_1$ to $\epsilon^d$.
\end{itemize}
\end{lemma}

\begin{proof}
(i)\ \
As we have seen in Example~\ref{eg:periodic},
every $\epsilon$-type periodic braid
is C-tight under the BKL Garside structure.
We now show that $\alpha^r$ is a BCMW-power of $\alpha$
using Theorem~\ref{thm:LL-CMW}.
Let $q=(n-1)/d$ and $k'=k/d$, then
$$
\INF(\alpha)=\INF(\epsilon^k)=\frac{n}{n-1}\cdot k
=\frac{nk'}q.
$$
Since $q$ is relatively prime to both $n$ and $k'$,
$nk'$ is relatively prime to $q$.
Because $\frac{n}{\gcd(nk',n)}=\frac{n}{n}=1$,
it suffices to show that $nk'\cdot r\equiv 1\bmod q$.
This follows from the following formula.
$$
nk'r=\frac{nkr}d=\frac{n(d-(n-1)s)}d=
\frac{n(d-dqs)}d=n(1-qs)=(1+dq)(1-qs)\equiv 1\bmod q.
$$

Since $\alpha$ is conjugate to $\epsilon^k$, the power
$\alpha^r$ is conjugate to
$
\epsilon^{kr}=\epsilon^{d-(n-1)s}
=(\epsilon^{n-1})^{-s} \epsilon^d
=\delta^{-ns}\epsilon^d$.

\smallskip\noindent
(ii)\ \
It is an easy consequence of (i).
We just use Proposition~\ref{prop:BCMW_power} that
$\gamma^{-1}\alpha\gamma=\epsilon^k$ if and only if
$\gamma^{-1}\alpha^r\gamma=\epsilon^{kr}$.
\end{proof}

\begin{corollary}
Let both $\alpha$ and $\beta$ be conjugate to $\epsilon^k$ for $k\ne 0$.
Let $d=\gcd(k,n-1)$, and $r$ and $s$ be integers such that $kr+(n-1)s=d$.
Let $\alpha_1=\delta^{ns}\alpha^r$ and $\beta_1=\delta^{ns}\beta^r$.
Then an $n$-braid $\gamma$ conjugates $\alpha$ to $\beta$
if and only if\/ it conjugates $\alpha_1$ to $\beta_1$.
In other words, the conjugacy search problem for $(\alpha,\beta)$
is equivalent to the conjugacy search problem for $(\alpha_1,\beta_1)$.
\end{corollary}

\begin{proof}
There is $\gamma_1\in B_n$ with $\beta=\gamma_1^{-1}\epsilon^k\gamma_1$,
which implies $\beta_1=\gamma_1^{-1}\epsilon^d\gamma_1$
by applying Lemma~\ref{lem:exp_red} to $(\beta, \beta_1)$.
Now apply Lemma~\ref{lem:exp_red} to $(\alpha, \alpha_1)$,
then we have
$$
\alpha=\gamma^{-1}\beta\gamma = \gamma^{-1}\gamma_1^{-1}\epsilon^k\gamma_1\gamma
\quad\mbox{if and only if}\quad
\alpha_1=\gamma^{-1}\gamma_1^{-1}\epsilon^d\gamma_1\gamma
=\gamma^{-1}\beta_1\gamma.
$$
\vskip-\baselinestretch\baselineskip
\end{proof}

From the above observation, to solve the conjugacy search
problem for periodic braids, it suffices to consider
the conjugacy classes of $\epsilon^d$ in $B_n$
only for the divisors $0<d< n-1$ of $n-1$
instead of $\epsilon^k$ for all integers $k$.
Before proving Proposition~\ref{prop:main},
we need to show a property of $\delta$.

\begin{lemma}\label{lem:delta}
Let $\alpha$ and $\beta$ be $n$-braids
such that $\delta=\alpha\beta=\beta\alpha$.
Then $\alpha$ is a power of $\delta$.
\end{lemma}

\begin{proof}
Since $\alpha\delta=\alpha(\beta\alpha)=(\alpha \beta)\alpha=\delta \alpha$,
the braid $\alpha$ belongs to the centralizer of $\delta$.
It is well-known that the centralizer of $\delta$ is the infinite cyclic group
generated by $\delta$ (see~\cite{BDM02} or~\cite{GW04}),
hence $\alpha$ is a power of $\delta$.
\end{proof}

\def\temp{
Recall that a simple element is a product of parallel descending cycles.
A descending cycle $[i_k, \ldots, i_2, i_1]$ in $B_n$ is defined originally
for the indices $i_1, \ldots, i_k$ with $1\le i_1 < i_2 < \cdots < i_k\le n$,
and indicates the positive word $a_{i_k i_{k-1}}\cdots a_{i_3 i_2} a_{i_2 i_1}$.
For notational simplicity, we will allow indices before making them congruent modulo $n$.
Namely, the form $[i_{j}+n,\ldots, i_1 + n, i_k, \ldots, i_{j+1} ]$
for any $j$ means the descending cycle $[i_k, \ldots, i_2, i_1]$.
For example, the form $[12, 11, 10, 9]$ in $B_{10}$ means
the descending cycle $[10,9,2,1]$.
}

For an $n$-braid $\alpha$,
let $\pi(\alpha)$ denote the induced permutation of $\alpha$.
We assume that the $n$-permutation $\pi(\alpha)$ acts
on $\{1,2,\ldots,n\}$ from right, and the expression
$k*\pi(\alpha)$ indicates the image of $k$ under
the action of  $\pi(\alpha)$.
For $0<d<n-1$,
notice that $\infs(\epsilon^d)= \inf(\epsilon^d)= d$ and
$\lens(\epsilon^d)= \len(\epsilon^d)= 1$ (by Lemma~\ref{lem:per-elt}),
and that the fixed point set of $\pi(\epsilon^d)$ is $\{1\}$,
hence every conjugate of $\epsilon^d$ has exactly one pure strand.
The next proposition shows how to find very efficiently
a conjugating element from $\epsilon^d$ to any given element
in the super summit set of $\epsilon^d$.

\begin{proposition}\label{prop:main}
Let $0< d<n-1$ be a divisor of\/ $n-1$, and let $q=(n-1)/d$.
Let $\alpha=\delta^d a$ be an $n$-braid conjugate to $\epsilon^d$,
having the $t$-th strand pure,
where $a \in\D\backslash\{ e, \delta\}$ and $1\le t\le n$.
Then the following hold.
\begin{enumerate}
\item[(i)] If the simple element $a$ has only one descending cycle, then
$\epsilon^d = \delta^{t-1}\alpha\delta^{-(t-1)}$.

\item[(ii)] If the simple element $a$ has more than one parallel descending cycles,
then at most $q-1$ iterations of partial cycling on a descending cycle of it
reduce the number of parallel descending cycles.
\end{enumerate}
\end{proposition}

\begin{proof}
(i)\ \
Suppose the simple element $a$ has only one descending cycle.
First, suppose $t=1$, that is, the first strand of $\alpha$ is pure.
Let $\pi(\delta^d)$ and $\pi(a)$ be the induced permutations of
$\delta^d$ and $a$.
Because $1*\pi(\delta^d a)=1$ and $1*\pi(\delta^d)=d+1$, we have
$$
1\stackrel{\pi(\delta^d)}{\longrightarrow}
d+1\stackrel{\pi(a)}{\longrightarrow}
1.
$$
Since $d+1\ne 1$, the number 1 is contained in the descending cycle of $a$.
On the other hand, because $\alpha=\delta^ua$ is conjugate to
$\epsilon^d=\delta^d [d+1,d,\ldots,1]$,
the exponent sum of $a$ is equal to that of $[d+1,d,\ldots,1]$.
This means that the descending cycle of the simple element $a$
has $d+1$ numbers.
Therefore $a$ is of the form
$$
[n_d,n_{d-1},\ldots,n_1,1],\qquad 1<n_1<\cdots<n_{d-1}<n_d \le n.
$$
Since $n_d*\pi(a)=1$, we obtain $n_d=d+1$ and hence $a=[d+1,d,\cdots,1]$.
This means that $\alpha=\epsilon^d$.

Now consider general cases.
If the $t$-th strand of $\alpha$ is pure,
then the first strand of $\delta^{t-1}\alpha\delta^{-(t-1)}$
is pure, hence, by the above argument,
$\epsilon^d =\delta^{t-1}\alpha\delta^{-(t-1)}$.

\medskip\noindent
(ii)\ \
Notice that $\alpha\in[\epsilon^d]^S$ because $\len(\alpha)=1$.
Suppose the simple element $a$ has more than one parallel descending cycles.
Among them, take any descending cycle, say $a_1$.
Then $a=a_1a_2$ for a non-identity simple element $a_2$.
We will show that at most $q-1$ iterations of partial cycling on $a_1$
reduce the number of parallel descending cycles.

Let $\psi=\tau^{d}$.
The partial cycling of $\alpha$ by $a_1$ is
$$
\alpha_1=\delta^{d}a_2\psi^{-1}(a_1).
$$
By Lemma~\ref{lem:e^d}, we know that
$[\epsilon^d]^S$ is closed under partial cycling.
This implies that $\alpha_1\in [\epsilon^d]^S$ and hence
$a_2\psi^{-1}(a_1)$ is a simple element.
If the number of parallel descending cycles is not changed, then
$$
a_2\psi^{-1}(a_1) = \psi^{-1}(a_1) a_2.
$$
Now we do partial cycling on $\alpha_1$ by $\psi^{-1}(a_1)$
and obtain
$$
\alpha_2=\delta^{d} a_2\psi^{-2}(a_1).
$$
By the same reason as above, $a_2\psi^{-2}(a_1)$ is a simple element.
If the number of parallel descending cycles is not changed, then
$$
a_2\psi^{-2}(a_1)= \psi^{-2}(a_1) a_2.
$$
Now assume that up to $q-1$ iterations of partial cycling on $a_1$ do not
decrease the number of parallel descending cycles of $a$.
Then
$\psi^{-k}(a_1)a_2 = a_2\psi^{-k}(a_1)$ for all $k=0,1,\ldots,q-1$,
which implies that
\begin{equation}\label{eq:comm}
\psi^j(a_2)\psi^i(a_1)=\psi^i(a_1)\psi^j(a_2),
\qquad \mbox{for all $i,j\in\{0,1,\ldots,q-1\}$}.
\end{equation}

We know that $\epsilon^d$ is P-minimal and C-tight (by Lemma~\ref{lem:e^d}),
hence $\alpha\in [\epsilon^d]^S=[\epsilon^d]^{St}$ (by Theorem~\ref{thm:pa-cy}).
Notice that $\INF(\epsilon^d)=nd/(n-1)=n/q$
and that $n$ and $q$ are relatively prime.
Therefore $\delta=\psi^{q-1}(a)\psi^{q-2}(a)\cdots\psi(a) a$
(by Lemma~\ref{lem:Primi_equiv}).
Combining with Equation~(\ref{eq:comm}), we have
$$
\delta=AB=BA
$$
where $A=\psi^{q-1}(a_1)\psi^{q-2}(a_1)\cdots\psi(a_1)a_1$
and $B=\psi^{q-1}(a_2)\psi^{q-2}(a_2)\cdots\psi(a_2)a_2$.
Since $a_1$ and $a_2$ are not the identity, we have $A,B\in\D\setminus\{e,\delta\}$.
In particular $A$ and $B$ cannot be a power of $\delta$,
which contradicts Lemma~\ref{lem:delta}.
\end{proof}

\section{Algorithms for the conjugacy problem for periodic braids}

Using the results in the previous sections,
we construct algorithms for solving the conjugacy problem
for periodic braids in $B_n^\BKLbr$.
The following is an overview of the algorithms we will describe
in this section.
\begin{itemize}
\item
Algorithms I and III are basic algorithms from which
the other algorithms are constructed.
Algorithm I provides an efficient method for powering
periodic elements in Garside groups,
and Algorithm III solves the CSP for
periodic $n$-braids conjugate to $\epsilon^d$,
where $d$ is a proper divisor of $n-1$.

\item
Algorithm II solves the CDP for periodic braids
and the CSP for $\delta$-type periodic braids.
Our solution to the CDP for periodic braids is
more efficient than Algorithm A of Birman,
Gebhardt and Gonz\'alez-Meneses~\cite{BGG06c}
because we use Algorithm I, the power conjugacy algorithm for periodic elements.
Our solution to the CSP for $\delta$-type periodic braids
is the same as Algorithm B of Birman,
Gebhardt and Gonz\'alez-Meneses~\cite{BGG06c}.

\item
Algorithm IV solves the CSP for $\epsilon$-type periodic braids.

\item
Algorithm V is a complete algorithm for the CDP and the CSP for
periodic braids.
\end{itemize}

Because Algorithm I works for periodic elements in
any Garside group, we describe it separately in \S5.1.
The other algorithms are described in \S5.2,
which work for the braid groups with the BKL Garside structure.
In \S5.3, we compare the complexities and
the necessary implementations
of our algorithms and those of Birman, Gebhardt and Gonz\'alez-Meneses in~\cite{BGG06c}.

Given a Garside group $G$, let $G^{+}$, $\Delta$ and $\D$ denote
the positive monoid, the Garside element
and the set of simple elements, respectively, of $G$.
Let $m$ be the smallest positive integer such that $\Delta^m$ is central in $G$.

\smallskip
Before going into algorithms,
let us first discuss how to represent elements of Garside groups
for inputs of algorithms.
The following two types of words are commonly used:
a word in the atoms and a word in the simple elements.
For example, an element in $B_n^\BKLbr$
can be represented by a word in the band generators
or by a word in the simple elements which
are products of parallel descending cycles.

In the following, we define words in the simple elements
in a little unusual way.
We explain the motivation briefly with an example.
Let $g$ be an element of $G$ represented by the word
$$
W=a_1^{k_1} \underbrace{\Delta\cdots\Delta}_u a_2^{k_2},
$$
where $u\ge 1$, $a_1,a_2\in\D$ and $k_1,k_2\in\{-1,1\}$.
Assume that $u$ is very large.
A natural algorithm for computing the normal form of $g$
would be as follows:
(i) collect $\Delta$'s in $W$ and obtain
$g=a_1^{k_1}\Delta^u a_2^{k_2}
=\Delta^u \tau^u(a_1)^{k_1}a_2^{k_2}$;
(ii) compute $\tau^u(a_1)$;
(iii) compute the normal form of $\tau^u(a_1)^{k_1}a_2^{k_2}$;
(iv) output the normal form of $g$ which is
the concatenation of $\Delta^u$ and the normal form of  $\tau^u(a_1^{k_1})a_2^{k_2}$.
Here, we remark the following two things.
First, the word length of $W$ is $u+2$, hence if we use the usual result
on the complexity for computing normal form, it will be $\mathcal O(u^2T)$
for some constant $T$, which is unnecessarily large.
Second, if $\Delta$'s are already collected so that $g$ is represented
by $a_1^{k_1}\Delta^u a_2^{k_2}$ then the complexity for computing
the normal form of $g$ is independent of $u$.
Therefore, in the following definition, we allow powers of
$\Delta$ to be contained in a word $W$ in the simple elements,
and we discard them when measuring the word length of $W$.

\begin{definition}
By a \emph{word in the simple elements} in a Garside group $G$, we mean the
following type of word $W$:
$$
W=\Delta^{u_0} a_1^{k_1} \Delta^{u_1} a_2^{k_2} \Delta^{u_2}
\cdots a_r^{k_r}\Delta^{u_r},\qquad
u_i\in\Z,~a_i\in\D,~k_i\in\{-1,1\}.
$$
Define $|W|_\simple=r$, the number of the simple elements $a_i$.
We use the notation $W^\atom$ (resp. $W^\simple$)
to indicate that $W$ is a word in the atoms
(resp. in the simple elements), when we want to make it more clear.
\end{definition}

For example, if a word $W$ is the normal form of an element $g$ in a Garside group,
then $|W|_\simple$ is the same as the canonical length of $g$.
Observe the following.
\begin{itemize}
\item
Every atom is a simple element.
Therefore, a word in the atoms can be regarded
as a word in the simple elements with the same word length.

\item
A word $W$ in the simple elements, $W=\Delta^{u_0} a_1^{k_1} \Delta^{u_1}
\cdots a_r^{k_r}\Delta^{u_r}$, can be transformed
to a word in the atoms by replacing each simple element with a product of atoms.
Let $V^\atom$ be such a transformed word.
Then we have the following inequalities:
$$
|W|_\simple
\le |V^\atom|
\le \Vert\Delta\Vert\cdot (|W|_\simple + \sum |u_i|),
$$
where $|V^\atom|$ denotes the word length of $V^\atom$.
The above formula shows that words in the simple elements
provide a more efficient way to implement
elements in Garside groups than words in the atoms.
\end{itemize}

For the algorithms in \S\ref{sec:Alg_power} and \S\ref{sec:Alg_BKL},
we assume that the elements of Garside groups are represented by
words in the simple elements, and we analyze their complexities
with respect to $|\cdot|_\simple$.

\subsection{Power conjugacy algorithm for periodic elements in Garside groups}
\label{sec:Alg_power}

In this subsection, we discuss complexities for algorithms in Garside groups,
and then give an efficient method for powering periodic elements.

We first recall the following notions in~\cite{Deh02}.
For simple elements $a$ and $b$, there is a
unique simple element $c$ such that $ac=a\vee_L b$.
Such an element $c$ is called the \emph{right complement} of $a$ in $b$.
Similarly, the \emph{left complement} of $a$ in $b$
is the unique simple element $c$ such that $ca=a\vee_R b$.
For a simple element $a$, let ${}^*a$ and $a^*$ denote
the left and right complements of $a$ in $\Delta$, respectively.
Therefore ${}^*a$ and $a^*$ are the unique simple elements satisfying
${}^*aa=\Delta=aa^*$, that is, ${}^*a=\Delta a^{-1}$ and $a^*=a^{-1}\Delta$.

\begin{definition}
Let $\Tlat$ be the maximal time for computing the following simple elements:
\begin{itemize}
\item $a\wedge_L b$ and $a\vee_L b$ from simple elements $a$ and $b$;
\item ${}^*a$ and $a^*$ from a simple element $a$;
\item $\tau^u(a)$ from a simple element $a$ and integer $0< u<m$.
\end{itemize}
\end{definition}

\begin{remark}
Because $\tau(a) = \Delta^{-1}a\Delta =(a^{-1}\Delta)^{-1}\Delta =(a^*)^*$,
we can compute $\tau(a)$ by computing the right complements twice.
Therefore, for $0<u<m$, we can compute $\tau^u(a)$
by computing right complements at most $2u$ times.
However, there are usually more efficient methods
for computing $\tau^u(a)$.

In $B_n^\Artinbr$, the simple elements are
in one-to-one correspondence with the $n$-permutations.
If $\theta$ is the $n$-permutation corresponding to a simple element $a$,
then the permutation corresponding to $\tau(a)$ is
$\theta'$ defined by $\theta'(i)= n+1-\theta(n+1-i)$ for $1\le i\le n$.
Moreover $\tau^2$ is the identity, hence for any integer $u$,
$\tau^u(a)$ can be computed in time $\mathcal O(n)$.
Note that the $a\wedge_L b$ for simple elements $a$ and $b$
can be computed in time $\mathcal O(n\log n)$~\cite{Thu92}.

In $B_n^\BKLbr$, the simple elements are products of parallel descending cycles.
If $[i_\ell,\ldots,i_1]$ is a descending cycle, then
$\tau^u([i_\ell,\ldots,i_1])=[i_\ell+u,\ldots,i_1+u]$, hence
for a simple element $a$, $\tau^u(a)$ can be computed in time
$\mathcal O(n)$.
\end{remark}

\begin{lemma}\label{lem:lt-op}
For simple elements $a$ and $b$ in a Garside group $G$,
the following operations can be done in time $\mathcal O(\Tlat)$.
\begin{itemize}
\item[(i)] $a\wedge_R b$ and $a\vee_R b$.
\item[(ii)] $\tau^u(a)$ for an integer $u$.
\item[(iii)] The left and the right complements of $a$ in $b$.
\item[(iv)] The normal form of\/ $ab$.
\end{itemize}
\end{lemma}

\begin{proof}
(i)\ \
It is known by~\cite[Lemma 2.5~(ii)]{Deh02} that
$$
a\wedge_R b={}^*(a^*\vee_L b^*)\quad\mbox{and}\quad
a\vee_R b =({}^*a\wedge_L {}^*b)^*.
$$
Therefore $a\wedge_R b$ can be computed by computing
two right complements, one lcm and then one left complement,
hence it can be computed in time $\mathcal O(\Tlat)$.
Similarly for $a\vee_R b$.

\smallskip\noindent
(ii)\ \
For any integer $u$, there is an integer $u_0$ such that
$u\equiv u_0\bmod m$ and $0\le u_0<m$.
Because $\tau^m$ is the identity, $\tau^u(a)=\tau^{u_0}(a)$.

\smallskip\noindent
(iii)\ \
Let $c$ be the right complement of $a$ in $b$, that is, $ac=a\vee_L b$.
Since
$$
{}^*(a\vee_L b) ac= {}^*(a\vee_L b) (a\vee_L b)=\Delta,
$$
we obtain that ${}^*(a\vee_L b) a$ is a simple element and
$$
c
=\bigl({}^*(a\vee_L b) a\bigr)^{-1}\Delta
=\bigl({}^*(a\vee_L b) a\bigr)^*.
$$
Therefore the element $c$ can be computed in time $\mathcal O(\Tlat)$.
Similarly for the left complement.

\smallskip\noindent
(iv)\ \
Let $a'b'$ be the normal form of $ab$. Then
$$
a'=\Delta\wedge_L (ab)=(aa^*)\wedge_L (ab)=a(a^*\wedge_L b)
\quad\mbox{and}\quad
b'=(a^*\wedge_L b)^{-1}b.
$$
It is obvious that $a'$ can be computed in time $\mathcal O(\Tlat)$.
Note that $b'$ is the right complement of $a^*\wedge_L b$ in $b$,
hence it can be computed in time $\mathcal O(\Tlat)$ by (iii).
\end{proof}

Recall that, for a positive element $g$ in $G$,
$\Vert g\Vert$ denotes the maximal word length of $g$ in the atoms in $G^{+}$.
The following lemma is well-known.
See~\cite{DP99} and \cite{BKL01}.

\begin{lemma}\label{lem:time}
Let $g$ be an element of a Garside group $G$.
\begin{itemize}
\item[(i)]
Let $g$ be given as a word $W$ in the simple elements with $l=|W|_\simple$.
Then the normal form of $g$ can be obtained in time $\mathcal O(l^2 \cdot\Tlat)$.

\item[(ii)]
Let $g$ be in the normal form.
Then the normal forms of the cycling and the decycling of\/ $g$
can be obtained in time $\mathcal O(\len(g) \cdot\Tlat)$.

\item[(iii)]
Let $g$ be in the normal form.
Then the total number of cyclings and decyclings
in order to obtain a super summit element is
$\mathcal O(\len(g) \cdot \Vert\Delta\Vert)$.
Therefore we can compute a pair $(g_1,h_1)$ such that
$g_1\in[g]^S$ is in its normal form and $h_1^{-1}g_1h_1=g$
in time $\mathcal O(\len(g)^2 \cdot\Vert\Delta\Vert \cdot\Tlat)$.
\end{itemize}
\end{lemma}

Note that in $B_n^\Artinbr$ one has $\Vert\Delta\Vert=\Vert\Delta_{(n)}\Vert=n(n-1)/2$ and
$\Tlat=n\log n$ by ~\cite{Thu92},
and that in $B_n^\BKLbr$ one has
$\Vert\Delta\Vert=\Vert\delta_{(n)}\Vert=n-1$ and $\Tlat=n$ by~\cite{BKL98}.

\medskip

\begin{definition}
If an element $g$ in a Garside group $G$ is periodic and belongs to
its super summit set, we call it a {\em periodic super summit element}.
\end{definition}

In order to solve the CDP/CSP for periodic elements in Garside groups,
we need an algorithm for powering a periodic element
and then computing a super summit element of that power.

\medskip
\underline{Power conjugacy algorithm}

\begin{itemize}
\item[]
INPUT: an integer $r\ge 1$ and a periodic super summit element $g$ in $G$.

\item[]
OUTPUT: a pair $(h,x)$ of elements in $G$ such that
$h\in[g^r]^S$ and $x^{-1}h x=g^r$.
\end{itemize}
A naive algorithm would be the following.
\begin{itemize}
\item[1.] Compute the normal form of $g^r$.
\item[2.] Apply iterated cycling and decycling to $g^r$
until a super summit element $h$ is obtained.
Let $x$ be the conjugating element obtained in this process
such that $x^{-1}hx=g^r$.
\item[3.] Return $(h,x)$.
\end{itemize}
Note that $\len(g)=\lens(g)\le 1$,
because $g$ is a periodic super summit element.
Because $\len(g^r)=r\len(g)=r$ in the worst case,
the complexity of the above algorithm
when $g$ is given in the normal form is
$$
\mathcal O(r^2\cdot\Tlat)+\mathcal O(r^2\cdot\Vert\Delta\Vert\cdot \Tlat)
=\mathcal O(r^2\cdot\Vert\Delta\Vert\cdot \Tlat).
$$

We will improve this algorithm so as to have complexity
$\mathcal O(\log r \cdot\Vert\Delta\Vert\cdot \Tlat)$.
Our idea is based on the repeated squaring algorithm in $\Z/n\Z$,
also known as exponentiation by squaring, square-and-multiply algorithm,
binary exponentiation or double-and-add algorithm.
We remark that our algorithm is interesting not only because
it gives an efficient method for powering,
but also because it exploits a recent result on
abelian subgroups of Garside groups~\cite{LL06c}.

Let us explain the repeated squaring algorithm in $\Z/n\Z$  briefly.
See~\cite[Page 8]{Coh93} or~\cite[Page 48]{Sho05} for more detail.
Let $a$, $n$ and $r$ be given large positive integers
from which we want to compute $a^r\bmod n$.
A naive algorithm for computing $a^r$ is to iteratively multiply by $a$
a total of $r$ times.
We can do better.
Let $r=k_0+2 k_1+2^2 k_2+\cdots+2^t k_t$ be the binary
expansion of $r$. Then we have the formula
$$
a^r=\prod_{k_i\ne 0} \left(a^{2^i}\right).
$$
The right hand side is a product of at most
$t+1=\lfloor \log_2 r\rfloor +1$ terms
and $a^{2^i}$ can be obtained by squaring $i$ times
as $a^{2^i}=(\cdots((a^2)^2)^2\cdots)^2$.
Therefore $a^r$ can be computed by $\mathcal O(\log r)$ multiplications.
This idea is implemented as follows.

\begin{algorithm}{Algorithm} (Repeated squaring algorithm in $\Z/n\Z$)\\
INPUT: positive integers $a$, $n$ and $r$.\\
OUTPUT: $a^r\bmod n$.

\begin{itemize}
\item[1.]
Compute the binary expansion
$r=k_0+2k_1+2^2k_2+\cdots+2^tk_t$, $k_i\in\{0,1\}$, of $r$.
\item[2.]
Set $x\leftarrow 1\bmod n$.
\item[3.]
For $i\leftarrow t$ down to 0, do the following.
\begin{itemize}
\item[] If $k_i=1$, set $x\leftarrow x^2a\bmod n$.
Otherwise, set $x\leftarrow x^2\bmod n$.
\end{itemize}
\item[4.] Return $x$.
\end{itemize}
\end{algorithm}

Let $r_i=\lfloor r/2^i\rfloor=k_i+2k_{i+1}+\cdots+2^{t-i}k_t$.
Then $r_i=2r_{i+1}+k_i$, hence
$a^{r_i}=(a^{r_{i+1}})^2a$ if $k_i=1$ and $a^{r_i}=(a^{r_{i+1}})^2$
if $k_i=0$.
At step 3 we are computing $a^{r_i}$
from $a^{r_{i+1}}$ for $i=t, t-1,\ldots,0$.

\begin{proposition}\label{prop:power}
Let $g$ be a periodic super summit element of a Garside group $G$,
being in the normal form.
For a positive integer $r$, there is an algorithm of complexity
$\mathcal O(\log r\cdot\Vert\Delta\Vert\cdot\Tlat)$
that computes a pair $(h,x)$ of elements of\/ $G$
such that $h\in[g^r]^S$ is in the normal form and $x^{-1}hx=g^r$.
\end{proposition}

\begin{proof}
Let $t=\lfloor \log_2 r\rfloor $, then the binary expansion of $r$ is as follows
$$
r=k_0+2 k_1+2^2 k_2+\cdots+2^t k_t,\qquad
\mbox{$k_i\in\{0,1\}$ for $i=0,\ldots,t$.}
$$
For $i\ge 0$, let
$$
r_i=\lfloor r/2^i\rfloor=k_i+2 k_{i+1}+2^2 k_{i+2}+\cdots+2^{t-i} k_t.
$$
Using reverse induction on $i$, we show that,
for each $i=t+1,t,\ldots,1,0$,
we can compute a triple $(g_i, h_i,x_i)$ such that
\begin{equation}\label{eq:ind}
g_i\in[g]^S, \qquad
h_i\in[g^{r_i}]^S, \qquad
x_i^{-1}g_ix_i=g, \qquad
x_i^{-1}h_ix_i=g^{r_i}.
\end{equation}
Notice that $(h_0,x_0)$ is the desired pair because $r_0=r$.
First, define
$$
g_{t+1}=g,\qquad
h_{t+1}=e,\qquad
x_{t+1}=e.
$$
It is obvious that $(g_{t+1},h_{t+1},x_{t+1})$
satisfies Equation~(\ref{eq:ind}) because $g^{r_{t+1}}=g^0=e$.
Assume that we have computed $(g_{i+1},h_{i+1},x_{i+1})$
for $0\le i\le t$.
Define
\begin{equation}\label{eq:len3}
h_i'=
\left\{\begin{array}{ll}
(h_{i+1})^2g_{i+1} & \mbox{if $k_i=1$,}\\
(h_{i+1})^2 & \mbox{otherwise.}
\end{array}\right.
\end{equation}
Note that $x_{i+1}$ conjugates $h_i'$ to $g^{r_i}$ because
$$x_{i+1}^{-1}h_i'x_{i+1} = \left\{
\begin{array}{ll}
(x_{i+1}^{-1}h_{i+1}x_{i+1})^2
(x_{i+1}^{-1}g_{i+1}x_{i+1})
= (g^{r_{i+1}})^2g
=g^{1+2r_{i+1}}
=g^{r_i} & \mbox{if } k_i=1,\\
(x_{i+1}^{-1}h_{i+1}x_{i+1})^2
= (g^{r_{i+1}})^2 =g^{2r_{i+1}}
=g^{r_i}
& \mbox{if } k_i=0.
\end{array}\right.$$
Apply iterated cycling and decycling to $h_i'$ until
a super summit element $h_i$ is obtained.
Let $y_i$ be the conjugating element obtained in this process such that
$h_i=y_ih_i'y_i^{-1}$.
Let
$$
x_i=y_ix_{i+1}
\quad\mbox{and}\quad
g_i=y_ig_{i+1}y_i^{-1}.
$$

Now we claim that $g_i$ is a super summit element.
Notice that $x_{i+1}$ conjugates $g_{i+1}$ and $h_i'$
to $g$ and $g^{r_i}$ respectively,
and that $g$ and $g^{r_i}$ commute with each other.
Therefore $g_{i+1}$ and $h_i'$ commute with each other,
that is, $g_{i+1}h_i'=h_i'g_{i+1}$.
In Lemma 3.2 of~\cite{LL06c}, the following is proved.
\begin{quote}
Let $g$ and $h$ be elements of a Garside group such that $gh=hg$.
Let $x$ be the conjugating element obtained in the process of applying
arbitrary iteration of cycling and decycling to $h$.
If $g$ is a super summit element, then so is $x^{-1}gx$.
\end{quote}
Therefore $g_i$ is a super summit element.
Note that $h_i$ is a super summit element by construction.
The element $x_i$ conjugates $g_i$ and $h_i$ to
$g$ and $g^{r_i}$ respectively, since
\begin{eqnarray*}
x_i^{-1}g_ix_i
&=& (y_ix_{i+1})^{-1}(y_ig_{i+1}y_i^{-1})(y_ix_{i+1})
= x_{i+1}^{-1}g_{i+1}x_{i+1}
= g,\\
x_i^{-1}h_ix_i
&=& (y_ix_{i+1})^{-1}(y_ih_i'y_i^{-1})(y_ix_{i+1})
= x_{i+1}^{-1}h_i'x_{i+1}
= g^{r_i}.
\end{eqnarray*}
Therefore $(g_i,h_i,x_i)$ satisfies Equation~(\ref{eq:ind}).

\medskip
Now we analyze the complexity of the above algorithm.
By definitions, both $h_{t+1} (=e)$ and $g_{t+1} (=g)$ are already in the normal form.
Assume that $h_{i+1}$ and $g_{i+1}$ are in the normal form for some $0\le i\le t$.

First we will show that one can compute the normal form of $h'_i$
from $(h_{i+1}, g_{i+1})$ in time $\mathcal O(\Tlat)$.
By the definition of $h'_i$ in (\ref{eq:len3}),
$h'_i$ is either $(h_{i+1})^2 g_{i+1}$ or $(h_{i+1})^2$.
Since both $h_{i+1}$ and $g_{i+1}$ are periodic super summit elements,
both $\len(h_{i+1})$ and $\len(g_{i+1})$ are at most 1.
Hence the number of non-$\Delta$ factors in the word representing $h'_i$
is at most 3, from which it follows that
we can compute the normal form of $h_i'$ in time $\mathcal O(\Tlat)$
and that $\len(h'_i)\le 3$.

Next we will show that one can compute the normal forms of $h_i$ and $g_i$
from $h'_i$ in time $\mathcal O(\Vert\Delta\Vert\cdot\Tlat)$.
Recall that $y_i$ is the conjugating element such that $h_i=y_ih_i'y_i^{-1}$
obtained in the process of applying iterated cycling and decycling to $h_i'$ until
a super summit element $h_i$ is obtained.
If $h'_i$ is already a super summit element, then $y_i=e$ hence we are done because
$h_i=h'_i$ and $g_i=g_{i+1}$.
Otherwise, $y_i\neq e$ and it is in fact given as a product
$y_i=y_{i,\ell} y_{i,\ell-1}\cdots y_{i,1}$
for some $\ell \ge 1$,
where each $y_{i,j}$ is a simple element or its inverse obtained in the process of \emph{each}
cycling or decycling from $h'_i$ to $h_i$.
Then $\ell\le\len(h'_i)\cdot\Vert\Delta\Vert\le 3\Vert\Delta\Vert$.
Using $y_{i,1}, y_{i,2}, \ldots, y_{i,\ell}$,
construct $h_{i,0},h_{i,1},\ldots,h_{i,\ell}$ and $g_{i,0},g_{i,1},\ldots,g_{i,\ell}$
recursively as
$$
h_{i, j+1} = y_{i,j+1}h_{i,j}y_{i,j+1}^{-1}\quad\mbox{and}\quad
g_{i, j+1} = y_{i,j+1}g_{i,j}y_{i,j+1}^{-1}\quad\mbox{for } 0\le j <\ell
$$
initializing $h_{i,0}=h'_i$ and $g_{i,0}=g_{i+1}$.
Then $h_{i,\ell}=h_i$ and $g_{i,\ell}=g_{i}$.
Notice that
$\len(h_{i,\ell})\le\len(h_{i,\ell-1})\le\cdots\le\len(h_{i,0})\le 3$
by the definition of $y_{i,j}$.
Notice also that
$\len(g_{i,\ell})=\len(g_{i,\ell-1})=\cdots =\len(g_{i,0})\le 1$
because $g_{i,j}h_{i,j}=h_{i,j}g_{i,j}$.

Since $h_{i,0}(=h'_i)$ is in the normal form and $y_{i,1}$ is a simple element or its inverse,
the number of non-$\Delta$ factors in the word representing
$h_{i,1}(=y_{i,1}h_{i,0}y_{i,1}^{-1})$ is at most 5.
Hence we can compute the normal form of $h_{i,1}$ from $h_{i,0}$
in time $\mathcal O(\Tlat)$.
In a recursive way,
we can compute the normal form of $h_{i,j+1}$ from $h_{i,j}$
in time $\mathcal O(\Tlat)$ for $j=0,\ldots,\ell-1$.
Since $\ell\le 3\Vert\Delta\Vert$,
we can compute the normal form of $h_{i}(=h_{i,\ell})$ from $h'_{i}(=h_{i,0})$
in time $\mathcal O(\Vert\Delta\Vert\cdot\Tlat)$.

The analogous proof works for computing the normal form of $g_{i,j+1}$ from
$(g_{i,j}, y_{i,j+1})$
in time $\mathcal O(\Tlat)$ for $j=0,\ldots,\ell-1$.
We just need to notice that the number of non-$\Delta$ factors in the word representing
$g_{i,j+1}(=y_{i,j+1}g_{i,j}y_{i,j+1}^{-1})$ is at most 3.
Thus we can compute the normal form of $g_{i}(=g_{i,\ell})$ from $g_{i+1}(=g_{i,0})$
in time $\mathcal O(\Vert\Delta\Vert\cdot\Tlat)$.

Therefore we can compute $(g_i,h_i,x_i)$ from $(g_{i+1},h_{i+1},x_{i+1})$
in time $\mathcal O(\Vert\Delta\Vert\cdot\Tlat)$ for $i$ from $t$ to $0$,
where $g_i$ and $h_i$ are in the normal form.
Since $t=\lfloor\log_2r\rfloor$, the whole complexity of the algorithm
is $\mathcal O(\log r\cdot \Vert\Delta\Vert\cdot\Tlat)$.
\end{proof}

The following is the algorithm discussed in Proposition~\ref{prop:power}.

\begin{algorithm}{Algorithm I}
(Power conjugacy algorithm for periodic elements in a Garside group $G$)\\
INPUT: a pair $(g,r)$ of a periodic super summit element $g\in G$
and a positive integer $r$, where $g$ is in the normal form.\\
OUTPUT: a pair $(h,x)$ of elements in $G$ such that $h\in[g^r]^S$ and $x^{-1}h x=g^r$,
where $h$ is in the normal form.

\begin{itemize}
\item[1.]
Compute the binary expansion
$r=k_0+2k_1+2^2k_2+\cdots+2^tk_t$, $k_i\in\{0,1\}$, of $r$.

\item[2.]
Set $g'\leftarrow g$, $h\leftarrow e$ and $x\leftarrow e$.

\item[3.]
For $i\leftarrow t$ down to 0, do the following.
\begin{itemize}
\item[3-1.]
If $k_i=1$, set $h'\leftarrow h^2g'$.
Otherwise, set $h'\leftarrow h^2$.
\item[3-2.]
Apply iterated cycling and decycling to $h'$ until
a super summit element $h$ is obtained.
Let $y$ be the conjugating element obtained in this process such
that $h=yh'y^{-1}$.
\item[3-3.]
If $\len(h)>1$, return ``\emph{$g$ is not a periodic element}''.
\item[3-4.]
Set $g'\leftarrow yg'y^{-1}$ and $x\leftarrow yx$.
\end{itemize}

\item[4.] Return $(h,x)$.
\end{itemize}
\end{algorithm}

\begin{remark}\label{rmk:Alg1}
Actually Algorithm I returns the desired pair $(h,x)$
if $\lens(g^k)\le 1$ for all $k$ even for a non-periodic element $g$.
In Algorithm I, Step 3.3 is used for the CDP when called by Algorithm II.
In any case, regardless of the summit length of $g^k$,
the complexity of Algorithm I is the same as the one in Proposition~\ref{prop:power}.
\end{remark}

\subsection{Algorithms in the braid groups with the BKL Garside structure}

Now we make an algorithm in $B_n^\BKLbr$ for solving  the CDP for periodic braids.
In~\cite{BGG06c}, Biman, Gebhardt and Gonz\'alez-Meneses proposed
the following algorithm.

\begin{algorithm}{Algorithm}
(Algorith A in~\cite{BGG06c} of Biman, Gebhardt and Gonz\'alez-Meneses)\\
INPUT: a word $W$ in the Artin generators representing an $n$-braid $\alpha$.\\
SUMMARY: determine whether $\alpha$ is periodic or not.
\begin{itemize}
\item[1.]
Compute the normal form of $\alpha^{n-1}$.\\
If it is equal to $\Delta^{2k}$,
return ``\emph{$\alpha$ is periodic and conjugate to $\epsilon^k$}''.

\item[2.]
Compute the normal form of $\alpha^n$.\\
If it is equal to $\Delta^{2k}$,
return ``\emph{$\alpha$ is periodic and conjugate to $\delta^k$}''.

\item[3.]
Return ``\emph{$\alpha$ is not periodic}''.
\end{itemize}
\end{algorithm}

\noindent
The word $W^{n-1}$ has word length $(n-1)l$ in the worst case,
where $l$ is the word length of $W$.
Therefore the complexity of the above algorithm is
$\mathcal O((ln)^2\cdot n\log n)=\mathcal O(l^2n^3\log n)$
as shown in~\cite[Proposition 5]{BGG06c}.
If one uses the BKL Garside structure in the above algorithm,
the complexity is reduced to $\mathcal O((ln)^2\cdot n)=\mathcal O(l^2n^3)$.
Using Algorithm~I and the fact that
every periodic braid has summit length at most 1,
we get a more efficient algorithm.

\begin{algorithm}{Algorithm II}
(Solving the CDP for periodic braids and the CSP for $\delta$-type periodic braids.)\\
INPUT: $\alpha\in B_n^\BKLbr$.\\
OUTPUT: ``\emph{$\alpha$ is not periodic}'' if $\alpha$ is not periodic;
``\emph{$\alpha$ is conjugate to $\epsilon^k$}''
if $\alpha$ is conjugate to $\epsilon^k$;
``\emph{$\alpha$ is conjugate to $\delta^k$ by $\gamma$}''
if $\gamma^{-1}\alpha\gamma=\delta^k$.
\begin{itemize}
\item[1.]
Compute the normal form of $\alpha$.

\item[2.]
Apply iterated cycling and decycling to $\alpha$ until
a super summit element $\beta$ is obtained.
Let $\gamma$ be the conjugating element obtained in this process such that
$\beta=\gamma^{-1}\alpha\gamma$.\\
If $\beta=\delta^k$,
return ``\emph{$\alpha$ is conjugate to $\delta^k$ by $\gamma$}''.\\
If $\len(\beta)>1$, return ``\emph{$\alpha$ is not periodic}''.

\item[3.]
Apply Algorithm I to $(\beta,n-1)$.\\
If it returns $\delta^{nk}$,
return ``\emph{$\alpha$ is conjugate to $\epsilon^k$}''.

\item[4.]
Return ``\emph{$\alpha$ is not periodic}''.
\end{itemize}
\end{algorithm}

\begin{theorem}\label{thm:CDP}
Let $\alpha$ be an $n$-braid represented by a word $W$ in the simple elements
of $B_n^\BKLbr$ with $|W|_\simple =l$.
Then there is an algorithm
of complexity $\mathcal O(l^2 n^2+n^2\log n)$
that decides whether $\alpha$ is periodic or not.
Further, if $\alpha$ is a $\delta$-type periodic braid,
then it decides that $\alpha$ is a $\delta$-type periodic braid
and computes a conjugating element $\gamma$ such that
$\gamma^{-1}\alpha\gamma=\delta^k$ in time $\mathcal O(l^2n^2)$.
\end{theorem}

\begin{proof}
Consider Algorithm II.
Step 1 computes the normal form of $\alpha$,
hence its complexity is $\mathcal O(l^2n)$ by Lemma~\ref{lem:time}~(i).
Step 2 computes a super summit element $\beta$,
hence its complexity is $\mathcal O(l^2n^2)$ by Lemma~\ref{lem:time}~(iii).
If $\alpha$ is a $\delta$-type periodic braid,
then Algorithm II stops here,
returning the conjugating element that conjugates $\alpha$ to $\delta^k$.
Therefore, at Step 3, we may assume that either
$\alpha$ is an $\epsilon$-type periodic braid
or it is not periodic.
In either case, Algorithm~I runs in time $\mathcal O(n^2 \log n)$ by Proposition~\ref{prop:power}
and Remark~\ref{rmk:Alg1}.
Therefore the total complexity of Algorithm~II is
$\mathcal O(l^2n^2+n^2\log n)$.
\end{proof}

Before considering the CSP for $\epsilon$-type periodic braids,
we study the case of periodic braids conjugate to $\epsilon^d$
for proper divisors $d$ of $n-1$.

Recall that a simple element in $B_n^\BKLbr$ is a product of parallel descending cycles.
A descending cycle $[i_k, \ldots, i_2, i_1]$ in $B_n$ is defined originally
for the indices $i_1, \ldots, i_k$ with $1\le i_1 < i_2 < \cdots < i_k\le n$,
and indicates the positive word $a_{i_k i_{k-1}}\cdots a_{i_3 i_2} a_{i_2 i_1}$.
For convenience, we will allow indices congruent modulo $n$.
Namely, the form $[i_{j}+n,\ldots, i_1 + n, i_k, \ldots, i_{j+1} ]$
for any $j$ means the descending cycle $[i_k, \ldots, i_2, i_1]$.
For example, the form $[12, 11, 10, 9]$ in $B_{10}$ means
the descending cycle $[10,9,2,1]$.

\begin{proposition}\label{prop:CSPe^d}
Let $\alpha$ be an $n$-braid in the normal form in $B_n^{\BKLbr}$.
If $\alpha\in[\epsilon^d]^S$ for a divisor $0<d< n-1$ of\/ $n-1$, then
there exists an algorithm of complexity $\mathcal O(n^2)$ that
computes a conjugating element $\gamma$
such that $\gamma^{-1}\alpha\gamma=\epsilon^d$.
\end{proposition}

\begin{proof}
Since $\lens(\epsilon^d)=1$ and $\alpha\in[\epsilon^d]^S$,
$\alpha=\delta^d a$ for some $a\in\D\backslash\{e,\delta\}$.
We will inductively construct sequences $\{\alpha_i\}_{i=0,\ldots,r}$
and $\{\gamma_i\}_{i=0,\ldots,r-1}$ of $n$-braids
for some $0\le r<d$ satisfying the following conditions.
\begin{itemize}
\item
$\alpha_i=\delta^d a_i\in [\epsilon^d]^S$
for a simple element $a_i\in\D\backslash\{e,\delta\}$ for each $i=0,\ldots,r$.

\item
$\gamma_i^{-1}\alpha_i\gamma_i=\alpha_{i+1}$ for $i=0,\ldots,r-1$.

\item
$a_0=a$ (and hence $\alpha_0=\alpha$).
The number of parallel descending cycles in $a_{i+1}$
is smaller than that in $a_i$ for $i=0,\ldots,r-1$.
The simple element $a_r$ has only one descending cycle.
\end{itemize}
Clearly, we can construct $\alpha_0$ by definition.
Suppose that we have constructed $\alpha_0,\ldots,\alpha_i$
and $\gamma_0,\ldots,\gamma_{i-1}$ for some $i$.
If $a_i$ has only one descending cycle,
then we already have constructed the desired sequences.
Therefore assume that $a_i$ has more than one parallel descending cycles.
Let $q=(n-1)/d$.
By Proposition~\ref{prop:main}~(ii),
at most $q-1$ iterations of partial cycling on
a descending cycle of $a_i$ reduce the number of parallel descending cycles in $a_i$.
Let $\alpha_{i+1}=\delta^d a_{i+1}$ denote the result
and let $\gamma_i$ be the conjugating element obtained in this process
such that $\gamma_i^{-1}\alpha_i\gamma_i=\alpha_{i+1}$.
Since $\alpha_i\in [\epsilon^d]^S$ and $[\epsilon^d]^S$ is closed under
partial cycling (by Lemma~\ref{lem:e^d}), $\alpha_{i+1}\in [\epsilon^d]^S$.
Notice that if one writes the simple element $a$ as a word in the band generators,
then the length is $d$, from which it follows that
there are at most $d$ parallel descending cycles in $a$.
Hence this process terminates in less than $d$ steps, that is, $r<d$.

Now we have the desired sequences, and $a_r$ has only one descending cycle.
By Proposition~\ref{prop:main}~(i), one has
$\alpha_r= \delta^{-(t-1)} \epsilon^d \delta^{t-1}$ for some $1\le t\le n$,
which means $a_r=[t+d,t+d-1,\ldots,t]$.
Let $\gamma=\gamma_0\gamma_1\cdots\gamma_{r-1}\delta^{1-t}$.
Then $\gamma^{-1}\alpha\gamma=\epsilon^d$.

\medskip
Because $r<d$ and we perform at most $q-1$ partial cyclings in order to obtain
$\alpha_{i+1}$ from $\alpha_i$ for $i=0,\ldots,r-1$,
the total number of partial cyclings in the whole process
is at most $d(q-1)<n$.
Because a partial cycling can be done in time $\mathcal O(n)$,
the complexity of this algorithm is $\mathcal O(n^2)$.
\end{proof}

The following is the algorithm discussed in Proposition~\ref{prop:CSPe^d}.

\begin{algorithm}{Algorithm III}
(Solving the CSP for $n$-braids conjugate to $\epsilon^d$.)\\
INPUT: the normal form $\delta^d a$ of an $n$-braid
$\alpha\in[\epsilon^d]^S$,
where $0<d<n-1$ is a divisor of $n-1$.\\
OUTPUT: an $n$-braid $\gamma$
such that $\gamma^{-1}\alpha\gamma=\epsilon^d$.

\begin{itemize}
\item[1.]
Set $\gamma \leftarrow e$.

\item[2.]
While $a$ has more than one parallel descending cycles, do the following.
\begin{itemize}
\item[2-1.]
Apply iterated partial cycling to $\alpha$ by a descending cycle of $a$
until we obtain a braid $\delta^d a'$
such that the number of parallel descending cycles in $a'$ is fewer than $a$.\\
Let $\gamma'$ be the conjugating element in this process such that
$\gamma'^{-1}(\delta^d a)\gamma' =\delta^d a'$.

\item[2-2.]
Set $a\leftarrow a'$ and  $\gamma\leftarrow \gamma\gamma'$.
\end{itemize}
\item[3.]
If $a$ has only one descending cycle,
say $[t+d,t+d-1,\ldots,t]$, set $\gamma\leftarrow \gamma \delta^{1-t}$.
\item[4.]
Return $\gamma$.
\end{itemize}
\end{algorithm}

\begin{figure}\tabcolsep=2pt
\begin{tabular}{*7c}
\includegraphics[scale=.55]{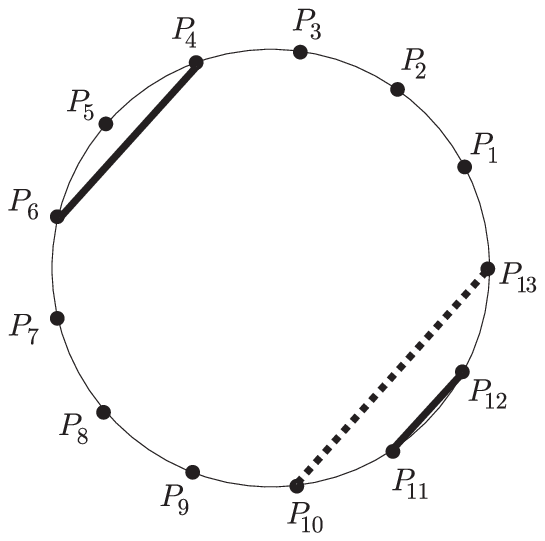} & \raisebox{13mm}{$\rightarrow$} &
\includegraphics[scale=.55]{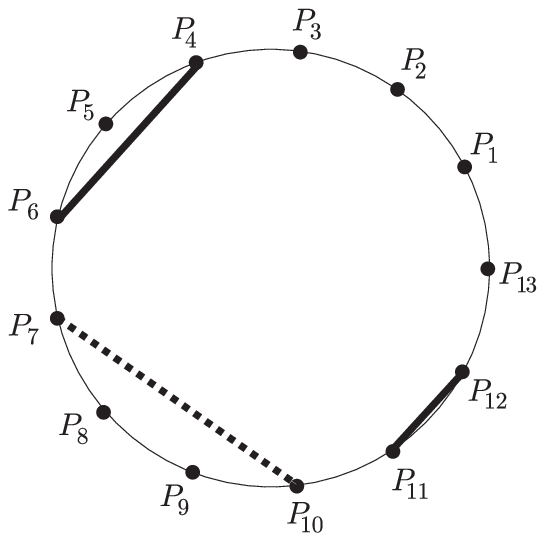} & \raisebox{13mm}{$\rightarrow$} &
\includegraphics[scale=.55]{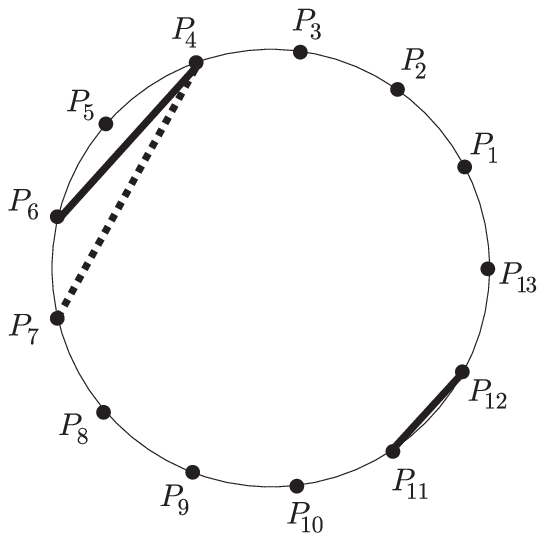} & \raisebox{13mm}{$=$} &
\includegraphics[scale=.55]{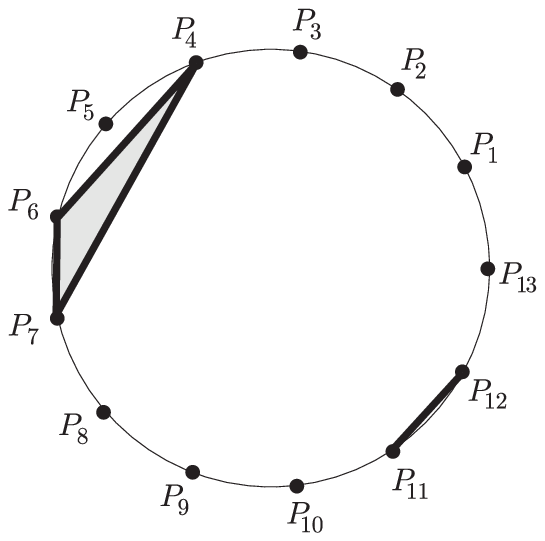}\\[5mm]
\includegraphics[scale=.55]{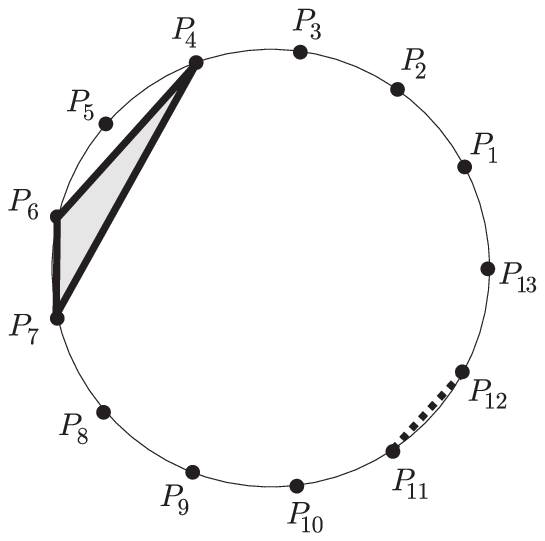} & \raisebox{13mm}{$\rightarrow$} &
\includegraphics[scale=.55]{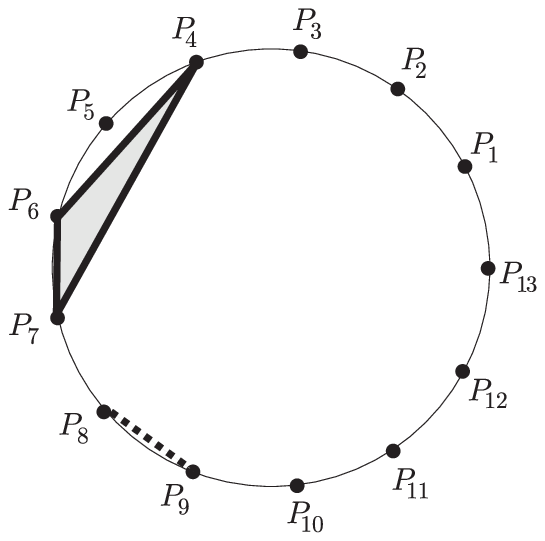} & \raisebox{13mm}{$\rightarrow$} &
\includegraphics[scale=.55]{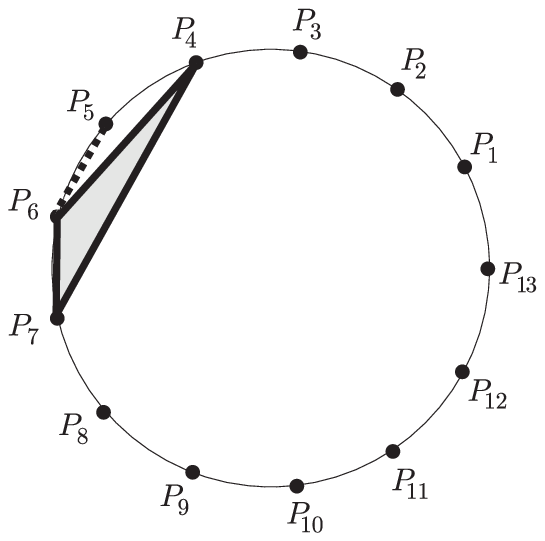} & \raisebox{13mm}{$=$} &
\includegraphics[scale=.55]{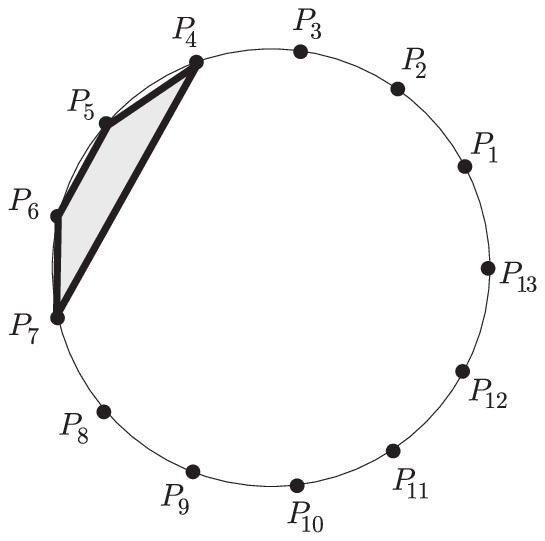}\\
\end{tabular}
\caption{The first row shows partial cyclings on the descending cycle
$[13,10]$ that is represented by dotted line.
After the partial cyclings $[13,10]\to [10,7]\to[7,4]$, the number
of parallel descending cycles is reduced by one.
Similarly, the second row shows partial cyclings
on the descending cycle $[12,11]$.
}\label{fig:move}
\end{figure}

\begin{example}\label{ex:largeUSS}
This example shows how Algorithm III transforms
an arbitrary periodic element
$\alpha\in [\epsilon^d]^S$ to $\epsilon^d$, where $0<d<n-1$ is a divisor of $n-1$.
See Figures~\ref{fig:e3conjB13} and~\ref{fig:move}.
Consider a 13-braid
$$
\alpha=\delta^3 [13,10][12,11][6,4].
$$
It is easy to see that $\alpha^4=\delta^{13} (=\epsilon^{12})$,
hence $\alpha$ is conjugate to $\epsilon^3=\delta^3[4,3,2,1]$.
Note that the simple element $[13,10][12,11][6,4]$ has three parallel descending cycles.

\begin{figure}
\includegraphics{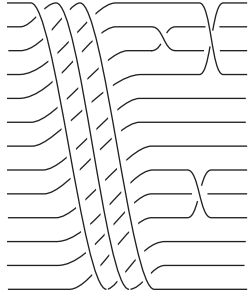}
\caption{The 13-braid $\alpha=\delta^3 [13,10][12,11][6,4]$}\label{fig:e3conjB13}
\end{figure}

\begin{itemize}
\item[(i)]
Iterate partial cycling on $[13,10]$ until it intersects another
descending cycle as follows:
$$
[13,10]\to [10,7]\to[7,4].
$$
Then the result is
$$
\alpha_1=\delta^3[12,11][6,4][7,4]=\delta^3[12,11][7,6,4].
$$
Note that $\alpha_1=b_1^{-1}\alpha b_1$, where $b_1=[10,7][7,4]=[10,7,4]$.

\item[(ii)]
Iterate partial cycling on $[12,11]$ until it intersects another
descending cycle as follows:
$$
[12,11]\to [9,8]\to[6,5].
$$
Then the result is
$$
\alpha_2=\delta^3[7,6,4][6,5]=\delta^3[7,6,5,4].
$$
Note that $\alpha_2=b_2^{-1}\alpha_1 b_2$, where $b_2=[9,8][6,5]$.

\item[(iii)]
Note that $\tau^{-3}(\alpha_2)=\delta^3[4,3,2,1]=\epsilon^3$.
Therefore $\beta^{-1}\alpha\beta=\epsilon^3$, where
\begin{eqnarray*}
\beta
&=&b_1b_2\delta^{-3}
=[10,7,4] [9,8][6,5]\delta^{-3}
=\delta^{-3}[7,4,1] [6,5][3,2].
\end{eqnarray*}
\end{itemize}
\end{example}

\begin{proposition}\label{prop:e^k}
Let $\alpha$ be an $n$-braid in the normal form in $B_n^{\BKLbr}$.
If $\alpha\in [\epsilon^k]^S$ for an integer $k$, then
there is an algorithm of complexity $\mathcal O(n^2\log n)$
that computes $\gamma$ such that $\gamma^{-1}\alpha\gamma=\epsilon^k$.
\end{proposition}

\begin{proof}
Let $u$ and $0\le v< n-1$ be integers satisfying $k=(n-1)u+v$,
and let $\alpha_0=\delta^{-nu}\alpha$.
Then $\alpha_0\in [\epsilon^v]^S$.
Since $\delta^n$ is central, for $\gamma\in B_n$
\begin{equation}\label{eq:n-1}
\gamma^{-1}\alpha\gamma=\epsilon^k
\quad\mbox{if and only if}\quad
\gamma^{-1}\alpha_0\gamma=\epsilon^v.
\end{equation}
If $v=0$, then $\gamma =e$, hence we may assume that $0< v< n-1$.

Compute $d=\gcd(v,n-1)$ and integers $r$ and $s$ such that
$$
d=vr+(n-1)s\quad\mbox{and}\quad 0\le r < n-1.
$$
Using Algorithm~I, compute $\gamma_1$ and the normal form of $\alpha_1$
such that $\alpha_1\in[\alpha_0^r]^S$ and
$$
\gamma_1^{-1}\alpha_1\gamma_1=\alpha_0^r.
$$
Let $\alpha_2=\delta^{ns}\alpha_1$, then $\alpha_2\in [\epsilon^d]^S$.
Apply Algorithm~III to $\alpha_2$ and obtain an $n$-braid $\gamma_2$ such that
$\gamma_2^{-1}\alpha_2\gamma_2=\epsilon^d$.
Lemma~\ref{lem:exp_red} shows that, for $\gamma\in B_n$,
\begin{equation}\label{eq:BCMW_power}
\gamma^{-1}\alpha_0\gamma=\epsilon^v
\quad\mbox{if and only if}\quad
\gamma^{-1}(\delta^{ns}\alpha_0^r)\gamma=\epsilon^d.
\end{equation}
Since $\delta^n$ is central,
$$
\epsilon^d = \gamma_2^{-1}\alpha_2\gamma_2
= \gamma_2^{-1}(\delta^{ns}\alpha_1)\gamma_2
= \gamma_2^{-1}\delta^{ns}(\gamma_1\alpha_0^r\gamma_1^{-1})\gamma_2
= \gamma_2^{-1}\gamma_1(\delta^{ns}\alpha_0^r)\gamma_1^{-1}\gamma_2.
$$
By Equations~(\ref{eq:n-1}) and~(\ref{eq:BCMW_power}),
$\gamma=\gamma_1^{-1}\gamma_2$ is the desired conjugating element.

\medskip
Now, let us analyze the complexity.
We can compute the integers $r$ and $s$ by using the extended Euclidean algorithm
which runs in time $\mathcal O(\log v\log n-1)=
\mathcal O((\log n)^2)$~\cite[Theorem 4.4 in page 60]{Sho05}.
By proposition~\ref{prop:power}, Algorithm~I with input $(\alpha_0,r)$ runs
in time $\mathcal O(n^2\log r)$.
Algorithm~III with $\alpha_2$ runs in time $\mathcal O(n^2)$ by
Proposition~\ref{prop:CSPe^d}.
Therefore the total complexity is
$$
\mathcal O((\log n)^2+n^2\log r+n^2)
=\mathcal O(n^2\log r)
=\mathcal O(n^2\log n).
$$
\vskip-\baselineskip
\end{proof}

The following is the algorithm discussed in Proposition~\ref{prop:e^k}.

\begin{algorithm}{Algorithm IV}
(Solving the CSP for $\epsilon$-type periodic braids)\\
INPUT: a pair $(\alpha,k)$, where $k$ is an integer and
$\alpha\in[\epsilon^k]^S$ is an $n$-braid in the normal form in $B_n^\BKLbr$. \\
OUTPUT: an $n$-braid $\gamma$ such that $\gamma^{-1}\alpha\gamma=\epsilon^k$.

\begin{itemize}
\item[1.]
Compute integers $u$ and $0\le v< n-1$ such that $k=(n-1)u+v$.

\item[2.]
Set $k \leftarrow v$ and
$\alpha\leftarrow \delta^{-nu}\alpha$

\item[3.]
If $k=0$, return $\gamma=e$.

\item[4.]
Compute $d=\gcd(k,n-1)$ and integers $0< r < n-1$ and $s$ such that $d=kr+(n-1)s$.

\item[5.]
Apply Algorithm I to $(\alpha, r)$.\\
Let $(\alpha_1,\gamma_1)$ be the output.
Then $\alpha_1\in[\alpha^r]^S$ is in the normal form,
and $\gamma_1^{-1}\alpha_1\gamma_1=\alpha^r$.\\
Set $\alpha_2\leftarrow\delta^{ns}\alpha_1$.
Then $\alpha_2$ belongs to the super summit set of $\epsilon^d$.

\item[6.]
Apply Algorithm III to $\alpha_2$.\\
Let $\gamma_2$ be the output, then $\gamma_2^{-1}\alpha_2\gamma_2=\epsilon^d$.

\item[7.]
Return $\gamma = \gamma_1^{-1}\gamma_2$.
\end{itemize}
\end{algorithm}

The following is the complete algorithm for the conjugacy problem
for periodic braids.

\begin{algorithm}{Algorithm V}
(The complete algorithm for the conjugacy problem for periodic braids)\\
INPUT: $\alpha\in B_n^\BKLbr$.\\
OUTPUT: ``\emph{$\alpha$ is not periodic}'' if $\alpha$ is not periodic;
``\emph{$\alpha$ is conjugate to $\epsilon^k$ by $\gamma$}''
if $\gamma^{-1}\alpha\gamma=\epsilon^k$;
``\emph{$\alpha$ is conjugate to $\delta^k$ by $\gamma$}''
if $\gamma^{-1}\alpha\gamma=\delta^k$.

\begin{itemize}
\item[1.]
Compute the normal form of $\alpha$.

\item[2.]
Apply iterated cycling and decycling to $\alpha$ until a super summit
element $\alpha_1$ is obtained.\\
Let $\gamma_0$ be the conjugating element in this process such that
$\gamma_0^{-1}\alpha\gamma_0=\alpha_1$.\\
If $\len(\alpha_1)>1$, return ``\emph{$\alpha$ is not periodic}''.

\item[3.]
Apply Algorithm II to $\alpha_1$.
\begin{itemize}
\item[3-1.]
If $\alpha_1$ is not periodic, return ``\emph{$\alpha$ is not periodic}''.

\item[3-2.]
If $\alpha_1$ is conjugate to $\delta^k$, Algorithm II gives
an element $\gamma_1$ such that $\gamma_1^{-1}\alpha_1\gamma_1=\delta^k$.\\
Set $\gamma\leftarrow\gamma_0\gamma_1$.\\
Return ``\emph{$\alpha$ is conjugate to $\delta^k$ by $\gamma$}''.

\item[3-3.]
If $\alpha_1$ is conjugate to $\epsilon^k$ for some $k$,
Algorithm II gives the exponent $k$.
\end{itemize}

\item[4.]
Apply Algorithm IV to $(\alpha_1,k)$. Let $\gamma_2$ be its output.\\
Set $\gamma\leftarrow \gamma_0\gamma_2$.\\
Return ``\emph{$\alpha$ is conjugate to $\epsilon^k$ by $\gamma$}''.
\end{itemize}
\end{algorithm}

\begin{proposition}\label{prop:CSP}
Let $\alpha$ be an $n$-braid given as a word $W$ in the simple elements
of $B_n^\BKLbr$ with $|W|_\simple =l$.
Then there is an algorithm of complexity $\mathcal O(l^2 n^2+n^2\log n)$
that decides whether $\alpha$ is periodic or not
and, if periodic, computes $\gamma\in B_n^\BKLbr$ such that
$\gamma^{-1}\alpha\gamma=\delta^k$ or
$\gamma^{-1}\alpha\gamma=\epsilon^k$.
\end{proposition}

\begin{proof}
It is not difficult to see that Algorithm V is the desired algorithm.
Now, let us analyze the complexity.
Step 1 computes the normal form,
hence its complexity is $\mathcal O(l^2n)$ by Lemma~\ref{lem:time}~(i).
Step 2 computes a super summit element,
hence its complexity is $\mathcal O(l^2n^2)$ by Lemma~\ref{lem:time}~(iii).
Step 3 applies Algorithm II to $\alpha_1$ whose word length is $\le 1$.
Therefore, its complexity is $\mathcal O(n^2\log n)$ by Theorem~\ref{thm:CDP}.
Step 4 uses Algorithm IV, hence the complexity is $\mathcal O(n^2\log n)$
by Proposition~\ref{prop:e^k}.
Therefore the total complexity is
$$
\mathcal O(l^2n+l^2n^2+n^2\log n)=
\mathcal O(l^2n^2+n^2\log n).
$$
\vskip-\baselinestretch\baselineskip
\end{proof}

\subsection{Remarks on efficiency of algorithms}

Here we compare our algorithms with the algorithms
of Birman, Gebhardt and Gonz\'alez-Meneses in~\cite{BGG06c}.
See Table~\ref{tab:BGG-alg} for their complexities,
the form of input words and necessary implementations.
Notice that the complexity of Algorithm~IV in Table~\ref{tab:BGG-alg}
is different from the one given in Proposition~\ref{prop:e^k}.
This is for the case where an input braid is not a super summit element and
not in the normal form like in Algorithm~C.

\begin{table}
{\small\tabcolsep=2pt
\def\arraystretch{1.2}
\begin{tabular}{|c||c|c|c|c|}
\multicolumn{5}{c}{(a) Our algorithms ($l = |W|_{\simple}$)}\\\hline
Problems & Algorithms & Complexity & Input word
& Necessary implementations \\\hline\hline
CDP & Algorithm II & $\mathcal O(l^2n^2+n^2\log n)$ & &\\\cline{1-3}
CSP for $\delta$-type & Algorithm II & $\mathcal O(l^2n^2)$
  & $W^\simple$ & Garside structure $B_n^\BKLbr$ \\\cline{1-3}
CSP for $\epsilon$-type & Algorithm IV
  & $\mathcal O(l^2n^2+n^2\log n)$ & & \\\hline
\multicolumn{5}{c}{}\\
\multicolumn{5}{c}{(b) Algorithms of Birman, Gebhardt
and Gonz\'alez-Meneses in~\cite{BGG06c} ($l = |W^\atom|$)}\\\hline
Problems & Algorithms & Complexity & Input word
& Necessary implementations \\\hline\hline
CDP & Algorithm A & $\mathcal O(l^2n^3\log n)$
  & & Garside structure $B_n^\Artinbr$\\\cline{1-3}\cline{5-5}
CSP for $\delta$-type & Algorithm B & $\mathcal O(l^3n^2)$
  & $W^\atom$ & Garside structures $B_n^\Artinbr$ \&\ $B_n^\BKLbr$ \\\cline{1-3}\cline{5-5}
CSP for $\epsilon$-type & Algorithm C & $\mathcal O(l^3n^2)$
  & & \begin{tabular}{c} Garside structures $B_n^\Artinbr$ \&\ $B_n^\BKLbr$\\[-.4em]
  bijections $P_{n,2}\leftrightarrows Sym_{2n-2}$\end{tabular}\\\hline
\end{tabular}}

\bigskip
\caption{Comparison of our algorithms with those in~\cite{BGG06c}}
\label{tab:BGG-alg}
\end{table}

In the paper~\cite{BGG06c}, Birman et.{} al.{}
proposed three algorithms for the conjugacy problem for periodic braids:
Algorithm A solves the CDP for periodic braids;
Algorithms B and C solve the CSP for $\delta$-type
and $\epsilon$-type periodic braids, respectively.
In this paper,
Algorithm II solves the CDP for periodic braids
and the CSP for $\delta$-type periodic braids,
and Algorithm IV solves the CSP for $\epsilon$-type periodic braids.

The main difference between the solutions of Biman et.{} al.{}
and ours is the way to solve the CSP for $\epsilon$-type periodic braids.
Algorithm C of Birman et.{} al.{} needs implementations
for the bijections between $P_{n,2}$ and $Sym_{2n-2}$, where
$P_{n,2}$ is a subgroup of $B_n$ consisting of all 2-pure braids, that is,
the $n$-braids whose induced permutations fix 2, and
$Sym_{2n-2}$ is the centralizer of $\delta_{(2n-2)}^{n-1}$ in $B_{2n-2}$.
It is known that both $P_{n,2}$ and $Sym_{2n-2}$ are isomorphic
to the Artin group of type $\B_{n-1}$, hence there exist
bijections from $P_{n,2}$ to $Sim_{2n-2}$ and vice versa.
Birman et.{} al.{} constructed the bijections explicitly
in~\cite{BGG06c}.

From Table~\ref{tab:BGG-alg}, our algorithms have
the following advantages.

\begin{itemize}
\item
The inputs of our algorithms are given as
words in the simple elements, while those of Birman et.{} al.{} are given as words in the atoms.
As we discussed at the beginning of this section,
it is more natural and more efficient
to represent elements in Garside groups as words in the simple elements.
For example, in the experiment in \S4 of~\cite{BGG06c},
Birman et.{} al.{} generate  a word in the simple
elements, not a word in the atoms.

\item
The complexities of our algorithms are lower than those of Birman et.{} al.{}
In the complexity of Algorithm C in Table~\ref{tab:BGG-alg}, $l \ge n$
unless the input braid is the identity element.

\item
Our algorithms require only implementations for Garside structure,
while the algorithms of Birman et.\ al.\ additionally require
implementations of the bijections between $P_{n,2}$ and $Sym_{2n-2}$.
\end{itemize}

We remark that Algorithms A and B of Birman el.{} al.{} can be revised as follows.
\begin{itemize}
\item
We can allow words in the simple elements
as the inputs of Algorithms A and B
without changing the complexity.
Let Algorithms A$'$ and B$'$ be the ones revised in this way.

\item
If we use the BKL Garside structure in Algorithm A$'$,
then the complexity is reduced from $\mathcal O(l^2n^3\log n)$
to $\mathcal O(l^2n^3)$.

\item
The complexity of Algorithm B$'$ is $\mathcal O(l^2n^2)$.
Algorithm B$'$ is the same as Algorithm II.
When analyzing the complexity of Algorithm B in~\cite{BGG06c},
they used that the time for computing the normal form after cycling or decycling
an element with canonical length $l$ is $\mathcal O(l^2n)$,
however it can be done in time $\mathcal O(ln)$.
\end{itemize}

\begin{center}
\small\bigskip
\begin{tabular}{|c||c|c|c|c|}\hline
Problems & Algorithms & Complexity & Input word
& Necessary implementations \\\hline\hline
CDP & Algorithm A$'$ & $\mathcal O(l^2n^3)$
  & \vbox to 0pt{\vss\hbox{$W^\simple $}\vskip -.8em}
  & \vbox to 0pt{\vss\hbox{Garside structure $B_n^\BKLbr$}\vskip -.8em}\\\cline{1-3}
CSP for $\delta$-type & Algorithm B$'$ & $\mathcal O(l^2n^2)$
  &  &  \\\hline
\end{tabular}
\bigskip
\end{center}

However, because of the transformations between $P_{n,2}$ and $Sym_{2n-2}$,
Algorithm C does not allow words in the simple elements as input without increasing
the complexity.

\end{document}